\documentclass{amsart}

\usepackage{amsmath,amssymb,amsthm}

\usepackage{graphicx,tikz,}

\newtheorem{theorem}{Theorem}

\newtheorem{proposition}{Proposition}

\newtheorem*{definition}{Definition}

\theoremstyle{remark}

\DeclareMathOperator{\supp}{supp}
\DeclareMathOperator{\diam}{diam}

\begin{document}

\title[]{Sums of Distances on Graphs and \\Embeddings into Euclidean Space}

\author[]{Stefan Steinerberger}
\address{Department of Mathematics, University of Washington, Seattle, WA 98195, USA}

\thanks{S.S. is supported by the NSF (DMS-2123224) and the Alfred P. Sloan Foundation.}

\begin{abstract} Let $G=(V,E)$ be a finite, connected graph. We consider a greedy selection of vertices:
given a list of vertices $x_1, \dots, x_k$, take $x_{k+1}$ to be any vertex maximizing the sum of distances to the
existing vertices and iterate: we keep adding the `most remote' vertex. The frequency with which the vertices of the graph appear in this sequence converges to a set of probability measures
with nice properties.
The support of these measures is, generically, given by a rather small number of vertices $m \ll |V|$. We prove that this suggests that the graph $G$ is
 at most '$m$-dimensional' by exhibiting an explicit $1-$Lipschitz embedding $\phi: G \rightarrow \ell^1(\mathbb{R}^m)$ with good properties.
 \end{abstract}

\maketitle

\vspace{-0pt}

\section{Introduction and Results}
\subsection{A Greedy Procedure.} Our original motivation was trying to understand the curious behavior of a simple procedure: given a finite, connected graph and given a list of vertices $x_1, \dots, x_{k}$ (note that this is a list: a vertex may appear multiple times), one
could try to extend the list by adding the vertex that is the furthest away in the sense of maximizing the sum over the distances to the existing vertices. This vertex may not be unique and we simply ask that $x_{k+1} \in V$ satisfies
$$ \sum_{i=1}^{k} d(x_i, x_{k+1}) = \max_{v \in V} \sum_{i=1}^{k} d(x_i, v).$$
The long-term behavior of this greedy procedure is interesting and connected to probability measures on $V$ with nice properties.
\begin{center}
\begin{figure}[h!]
\begin{tikzpicture}
\node at (0,0) {\includegraphics[width=0.25\textwidth]{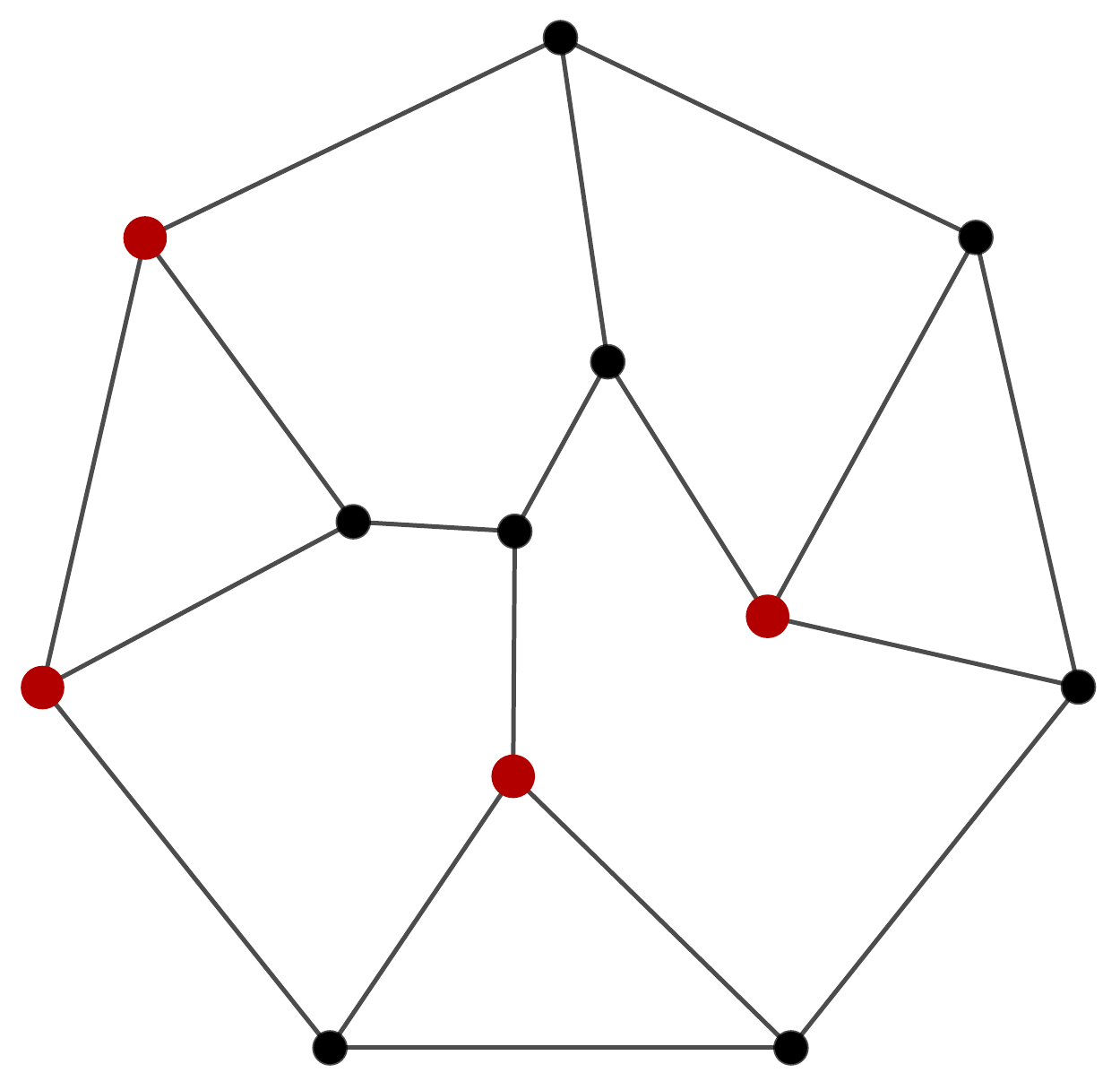}};
\node at (5,0) {\includegraphics[width=0.3\textwidth]{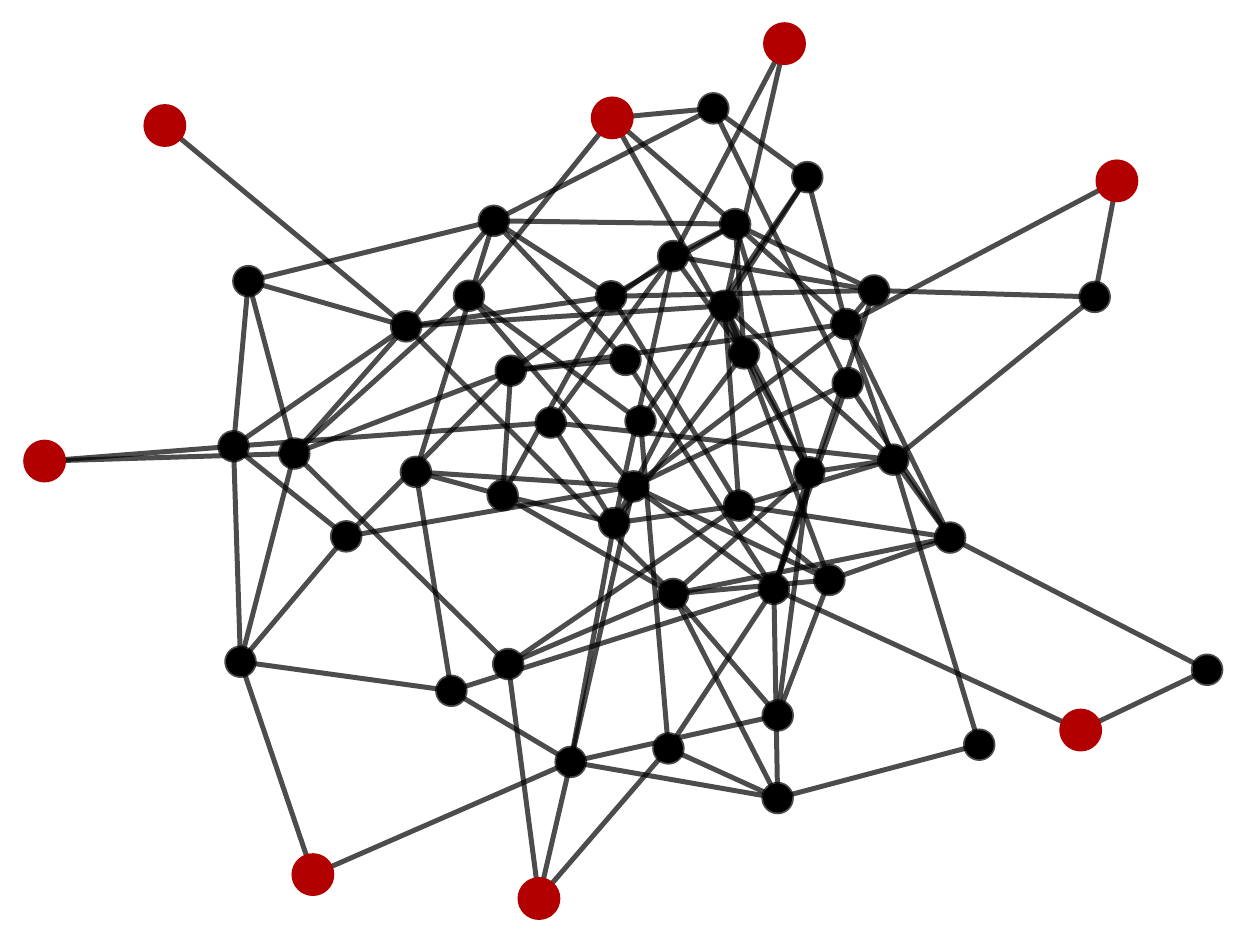}};
\end{tikzpicture}
\vspace{-10pt}
\caption{The Frucht graph (left) and an Erd\H{o}s-Renyi random graph on $n=50$ vertices. The rule ends up only selecting the red vertices (not necessarily with equal frequency).}
\end{figure}
\end{center}

Some experimentation suggests the frequency with which vertices arise
seems to quickly stabilize to a type of limiting distribution. We also observe that, typically, any such limiting measure $\mu$ seems to be supported on a much smaller subset of vertices. This is not always the case and it is possible, albeit fairly rare, that the limiting measure is actually given by the uniform measure on the vertices (see \S 2.1). However, for most graphs (both highly structured and random), the measure concentrates on a small subset.  One could think of this as a procedure favoring subsets of vertices that are at great distance from each other. The purpose of this paper is to demonstrate that the procedure has a number of interesting properties and applications (most notably inducing a graph embedding).

\begin{center}
\begin{figure}[h!]
\begin{tikzpicture}
\node at (0,0) {\includegraphics[width=0.4\textwidth]{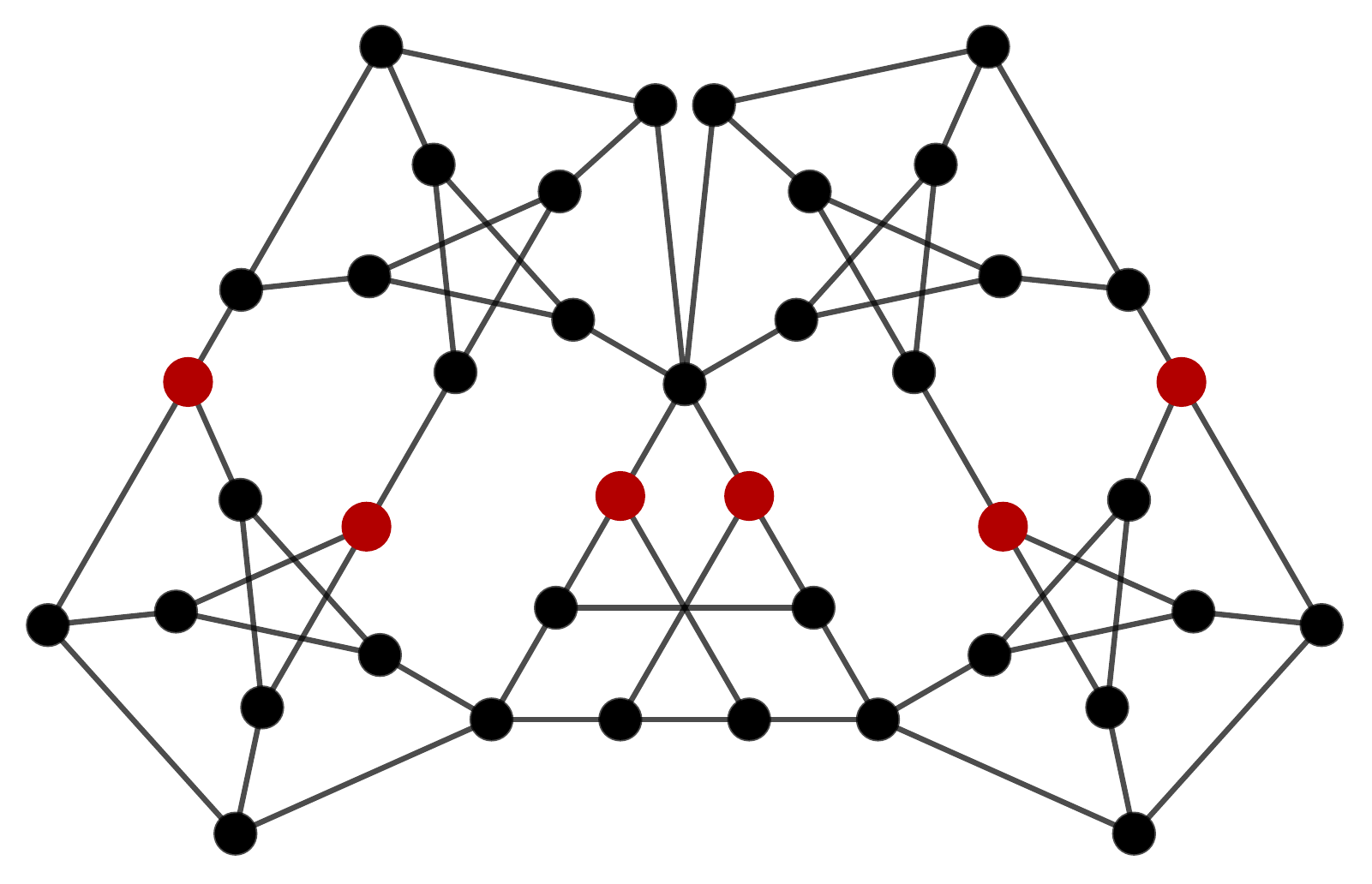}};
\node at (6,0) {\includegraphics[width=0.3\textwidth]{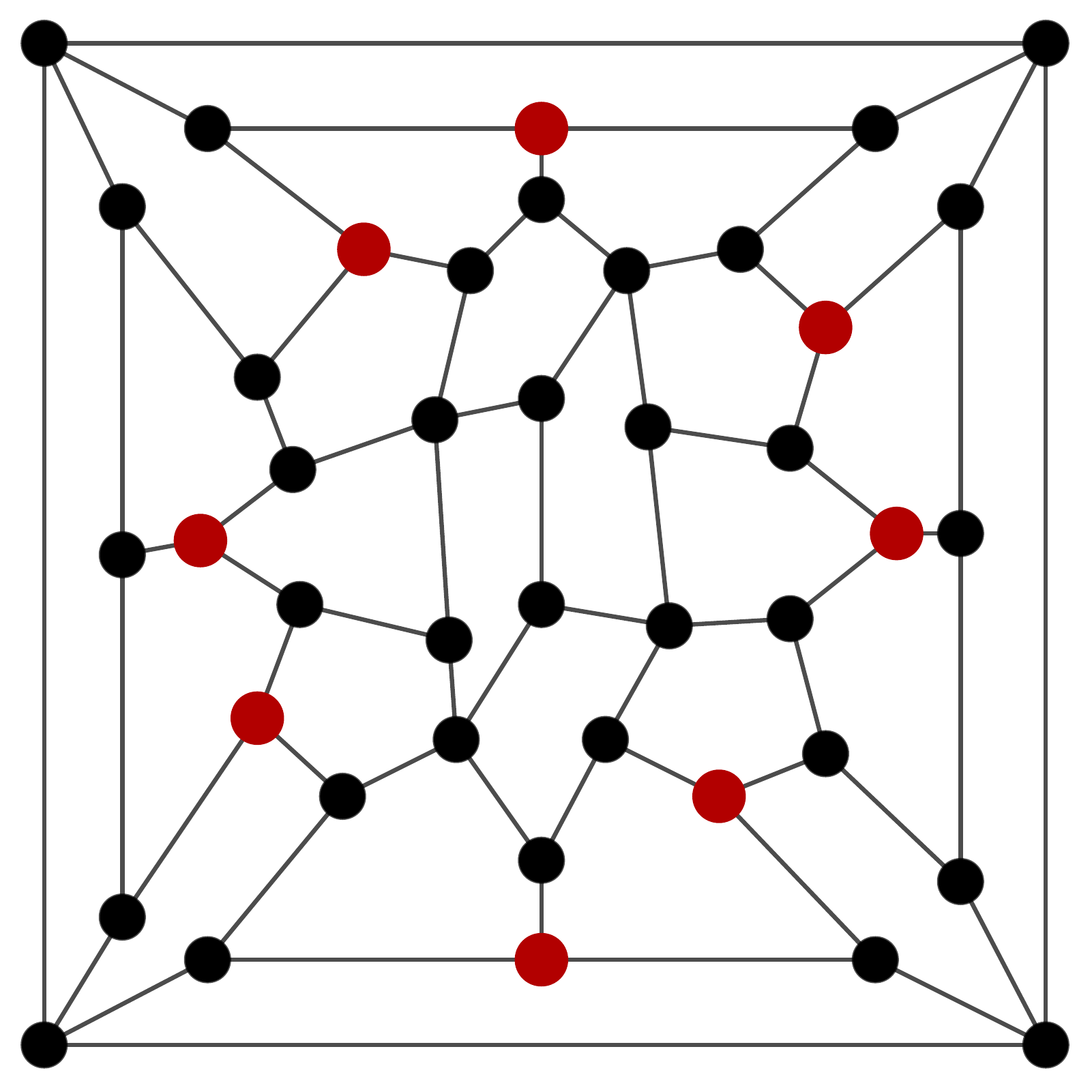}};
\end{tikzpicture}
\vspace{-5pt}
\caption{Thomassen 41 and Wiener-Araya.}
\end{figure}
\end{center}
\vspace{-15pt}

%\begin{center}
%\begin{figure}[h!]
%\begin{tikzpicture}
%\node at (0,0) {\includegraphics[width=0.4\textwidth]{ex3}};
%\node at (6,0) {\includegraphics[width=0.4\textwidth]{ex5}};
%\end{tikzpicture}
%\vspace{-15pt}
%\caption{Thomassen 60 (left) and Thomassen 105 (right).}
%\end{figure}
%\end{center}

\subsection{Long-time behavior.} Given an initial list $x_1, \dots, x_k \in V$, we use the rule 
$$ \sum_{i=1}^{k} d(x_{k+1}, x_i) = ~\max_{v \in V} \sum_{i=1}^{k} d(v, x_i).$$
to generate an infinite sequence of vertices. Note that an initial list does not necessarily specify a unique sequence but rather a family of sequences (since the maximum might be assumed in multiple vertices). We will not distinguish between different elements of the same family and always refer to them as `a sequence'. The empirical density of vertices of such a sequence will be shown to approach balanced probability distributions which we now define.
\begin{definition}
A probability measure $\mu$ on the vertices $V$ is \emph{balanced} if
$$ \mu(w) > 0 \implies \sum_{u \in V} d(w,u) \mu(u) = \max_{w \in W} \sum_{u \in V} d(w, u) \mu(u).$$
\end{definition}
One interpretation is as follows: for any probability measure $\mu$ on the vertices, we may introduce the (transport cost) function $T:V \rightarrow \mathbb{R}$ 
$$ T(w) = \sum_{u \in V} d(w,u) \mu(u).$$
$T(w)$ is the transport cost of sending all the mass of $\mu$ to the vertex $w$ under the assumption that transporting $\varepsilon > 0$ units of mass across one edge costs $\varepsilon$ transport cost -- this is also known as the Wasserstein $W^1$ cost or {Earth Mover Distance}. A probability measure $\mu$ on the vertices $V$ is said to be \textit{balanced} if it has the property that whenever $\mu(w) > 0$, then $T$ assumes a global maximum in $w$. Poetically put, a measure is balanced if the points where probability mass can actually be found are simultaneously the vertices for which global transport of the entire measure to a single vertex is the most expensive.
We can now describe the asymptotic behavior of the greedy vertex selection: it is assumed that we start with an arbitrary initial list of vertices and then turn the list into
an infinite sequence by picking vertices maximizing the sum of distances.

%Identifying probability measures $\mu$ on $V$ with probability vectors $v \in \mathbb{R}^n$, we have a natural metric on the space of probability measures by setting $d(\mu, \nu) = \| \mu - \nu\|_{\ell^{\infty}}$. Since the space is finite-dimensional, the choice of norm does not matter (up to constant).

\begin{theorem}[Long-time behavior]
Let $(x_k)_{k=1}^{\infty}$ denote an infinite sequence of vertices obtained by the greedy procedure and let $\mu_m$ denote the empirical probability measure of the first $m$ vertices. Then
\begin{enumerate}
\item there exists $\diam(G)/2 \leq \alpha \leq \diam(G)$ such that
$$ \lim_{m \rightarrow \infty} \sum_{v,w \in V} \mu_m(v) d(v,w) \mu_m(w) \rightarrow \alpha$$
\item  the maximal transport cost converges to $\alpha$
$$ \lim_{m \rightarrow \infty} \max_{v \in V} \sum_{w \in V} d(v, w) \mu_m(w) = \alpha $$
\item and for any $\varepsilon >0$ and all $m$ sufficiently large 
$$ \mu(v) > \varepsilon \quad \implies \quad  \sum_{w \in V} d(v, w) \mu_m(w) \geq \alpha - \varepsilon.$$
\end{enumerate}
\end{theorem}
This implies that any convergent subsequence of $\mu_m$ (which exist due to compactness) has to converge to a balanced measure.
One consequence is that at least one balanced probability measures always exists (though this is not difficult to show by other means, see Proposition 1).
Another consequence is that if one wishes to create a measure that is close to a balanced measure, one can simply use the greedy selection procedure
and is guaranteed to end up close to such a measure.

\subsection{Embedding $G=(V,E)$ into $\mathbb{R}^m$.}
If a balanced measure $\mu$ is supported on a small number $m \ll n$ of vertices, then this implies that the graph behaves approximately like an $m-$dimensional object in a way that we will make now precise. Our goal will be to find an embedding $\phi:V \rightarrow \mathbb{R}^m$ that is $1-$Lipschitz with respect to the graph distance: nearby points get mapped to nearby points. This is easy: one could simply map all vertices to the same point. We thus require an additional condition that the graph is not mapped into too small a region.

\begin{theorem}[Graph embedding]
Let $G=(V,E)$ be connected and let $\mu$ be a balanced measure on $m$ vertices. There exists $\phi:V \rightarrow \mathbb{R}^m$ such that
for all $u,v \in V$
$$ \| \phi(u) - \phi(v)\|_{\ell^1(\mathbb{R}^m)} \leq d(u,v).$$
Moreover, for some $\diam(G)/2 \leq \alpha \leq \diam(G)$, $\phi$ sends $\supp \mu$ to 
$$ \phi(\supp \mu) \subset \left\{x \in \mathbb{R}_{\geq 0}^m: x_1 + x_2 + \dots + x_m = \alpha\right\},$$
and points in $ \phi(\supp \mu)$
are, on average, not too close: for all $v \in \supp \mu$
$$ \frac{1}{\# \supp \mu} \sum_{ w \in \supp \mu} \| \phi(v) - \phi(w) \|_{\ell^{\infty}(\mathbb{R}^m)} \geq  \frac{\diam(G)}{2m}.$$
\end{theorem}
\newpage
\textbf{Remarks.}
\begin{enumerate}
\item The embedding $\phi$ is explicit: if $\mu$ is supported on $w_1, \dots, w_m \in V$, then
$$ \phi(v) = \left( \mu(w_1) d(w_1, v), ~ \dots, ~\mu(w_m) d(w_m, v) \right) \in \mathbb{R}^m.$$
This can be thought of as a triangulation by $\supp \mu$ weighted by $\mu$. If $\mu(w_j)$ is small,
then one can obtain another embedding of comparable quality in a lower dimension by omitting this
coordinate (see below for an example).\\
\item The inequality is sharp up to constants (see below for an example). Since the embedding is 1-Lipschitz, it would be interesting to obtain lower
bounds on the sums
$$ \sum_{ u,v \in \supp \mu} \| \phi(u) - \phi(v) \|_{\ell^{1}(\mathbb{R}^m)} \quad \mbox{and} \quad \sum_{u,v \in V} \| \phi(u)-\phi(v)\|_{\ell^1(\mathbb{R}^n)}$$
since the embedding is into $\ell^1(\mathbb{R}^m)$ rather than $\ell^{\infty}(\mathbb{R}^m)$. Is it possible to give conditions on $G=(V,E)$ under which $\diam_{\ell^1}(\phi(V)) \geq c \cdot \diam(G)$? An example shown below shows that this is not always the case.\\
\item There is no sense in which this $m$ would be minimal: it is certainly conceivable that the graph embedding $\phi(V)$ is
contained in a lower-dimensional subset of $\mathbb{R}^m$ and it is easy to construct examples for which this happens. The size of $\supp(\mu)$ provides an upper bound on the dimension. We also note that graphs can support balanced measures whose supports have a different cardinality (see also \S 2.3).\\
\end{enumerate}

\textbf{Graph Embedding: Examples.}\\

\textit{1. Combinatorial Example.} The Zamfirescu 75 graph \cite{zam} looks rather complicated (see Fig. 3). However, creating the sequence of vertices we quickly
find that the game ends up jumping between only three vertices. Computing the arising embedding into $\mathbb{R}^3$ demonstrates that there is a rather
simple underlying structure underlying the graph. Note that the embedding is not injective and does in fact collapse different vertices onto the same point in
$\mathbb{R}^3$.
\begin{center}
\begin{figure}[h!]
\begin{tikzpicture}
\node at (0,0) {\includegraphics[width=0.36\textwidth]{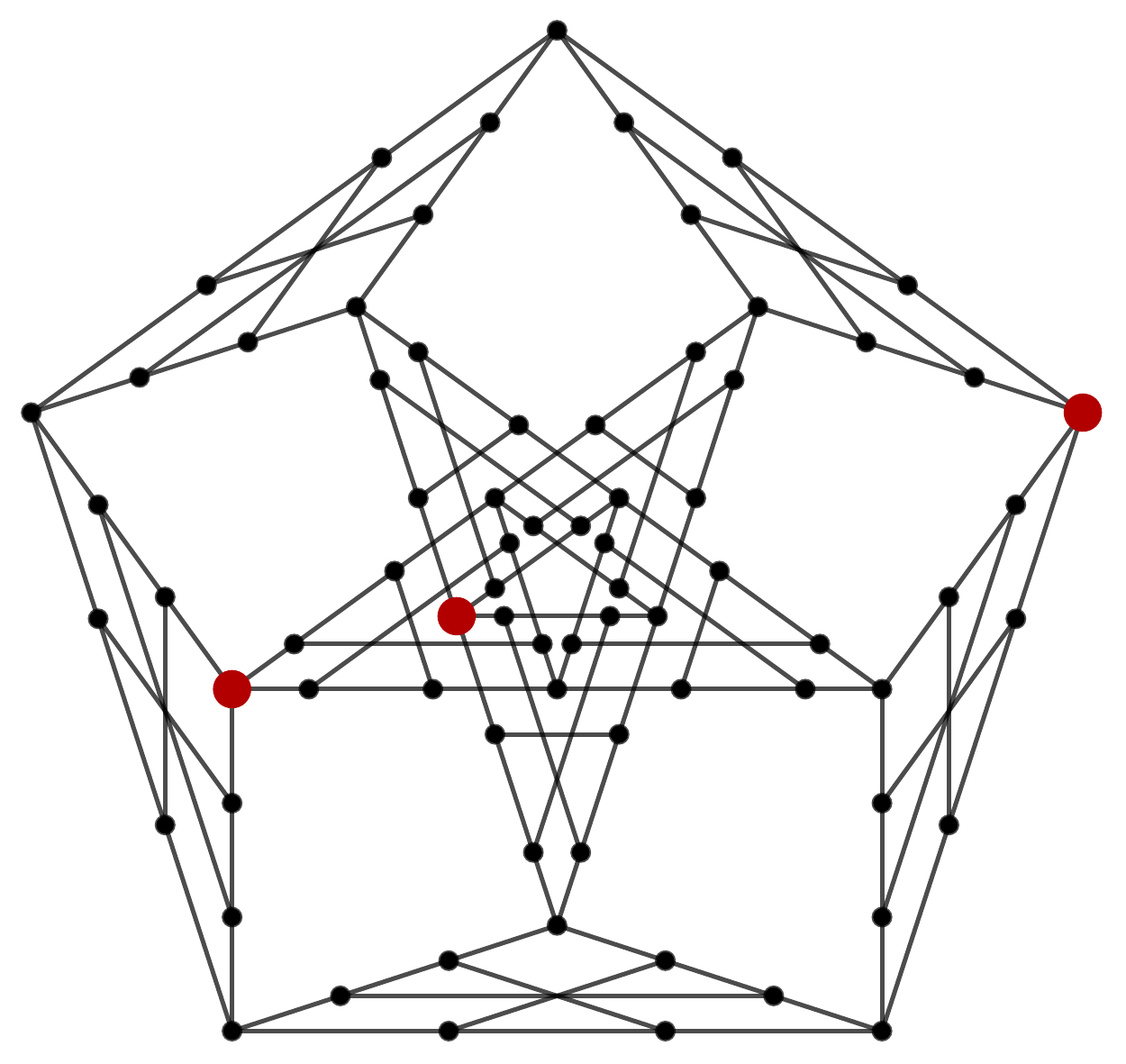}};
\node at (6,0) {\includegraphics[width=0.4\textwidth]{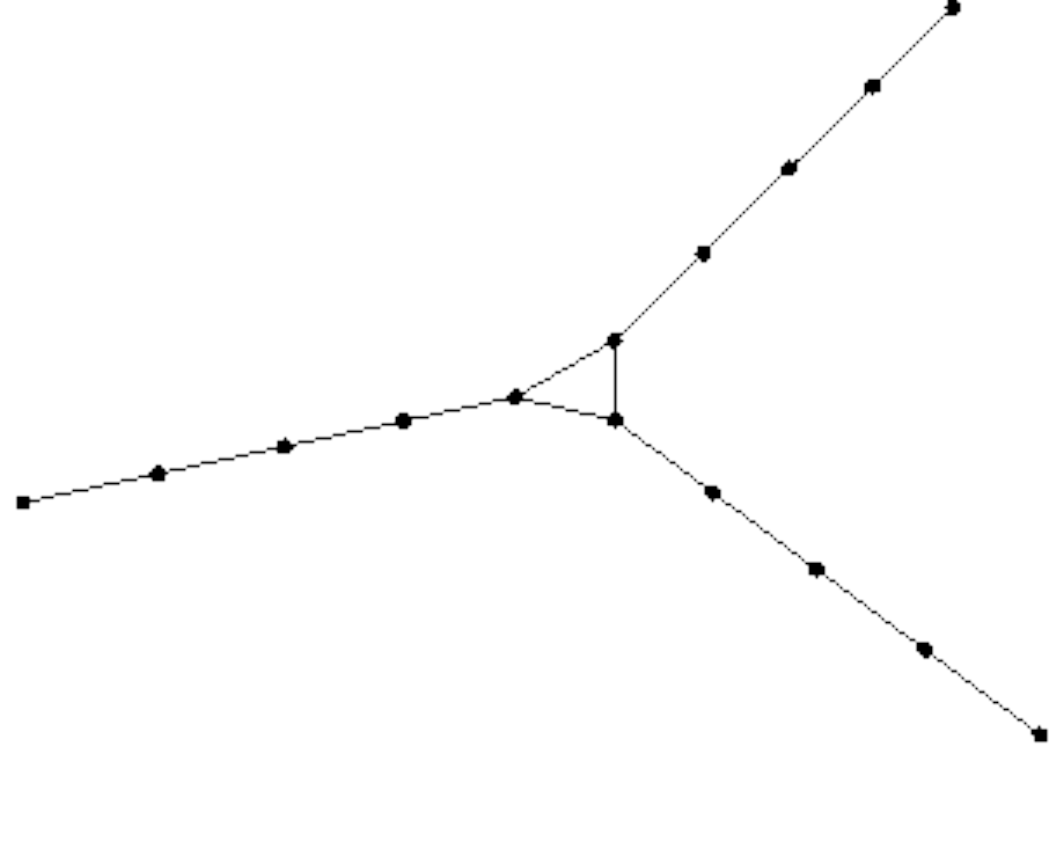}};
\filldraw[red] (3.6,-0.4) circle (0.06cm);
\filldraw[red] (8,1.95) circle (0.06cm);
\filldraw[red] (8.44,-1.5) circle (0.06cm);
\end{tikzpicture}
\vspace{-5pt}
\caption{Zamdirescu 75 (left) and the arising embedding (right, with lines connecting $\phi(v)$ and $\phi(w)$ for all $(v,w) \in E$).}
\end{figure}
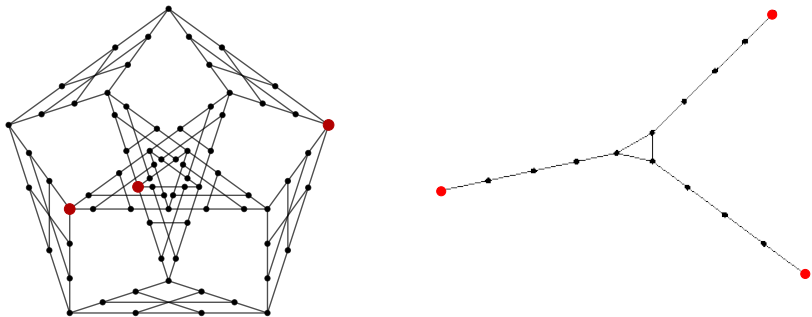
\end{center}
\vspace{-10pt}
This is an interesting example where the graph is actually combinatorially somewhat simpler than it may at first glance appear; this is
reflected in particularly simple long-term behavior of the game and a correspondingly simple embedding.\\

\textit{2. Gaussian Point Cloud.}  Another type of example is shown
in Fig. 4: these points are obtained from taking nearest-neighbor connections between Gaussian point clouds.
What we observe is that for these types of examples it seems that the embedding is fairly close to an isometry for a majority of vertices: 
up to a universal constant $c_G>0$ it seems that for the vast majority of pairs of vertices $(u,v) \in V$ (say $95\%$ of $V \times V$), we have
$$ \frac{1}{2} \leq c_G \frac{d(u,v)}{\| \phi(u) - \phi(v)\|_{\ell^1}} \leq \frac{3}{2}.$$
It would be interesting if such a statement could be made precise for, say, certain types of random graphs.
The purpose of this example is to illustrate that the Theorem does appear to actually lead to bi-Lipschitz embeddings with good constants even in somewhat rough settings.

\begin{center}
\begin{figure}[h!]
\begin{tikzpicture}
\node at (0,0) {\includegraphics[width=0.44\textwidth]{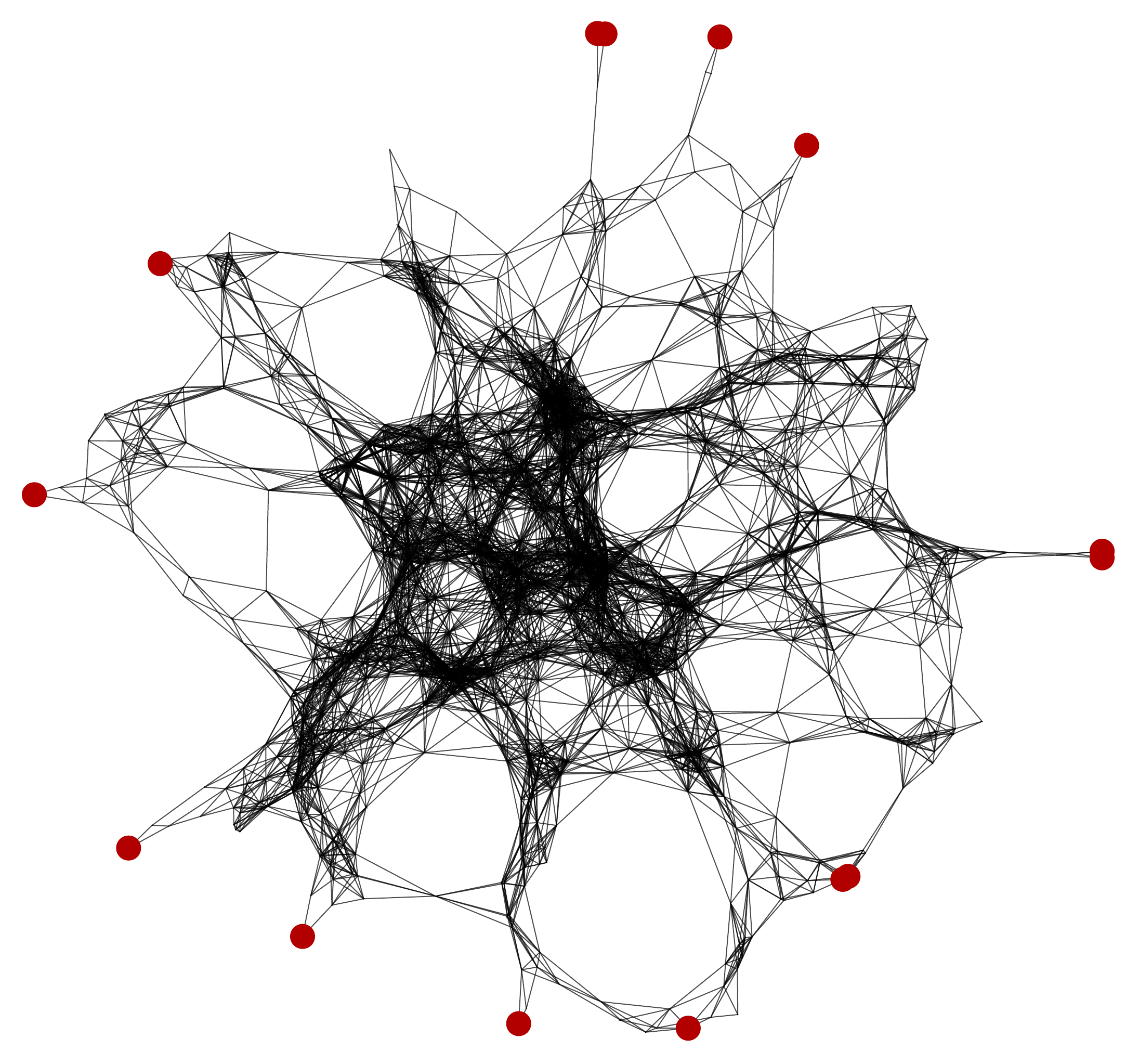}};
\node at (6,0) {\includegraphics[width=0.4\textwidth]{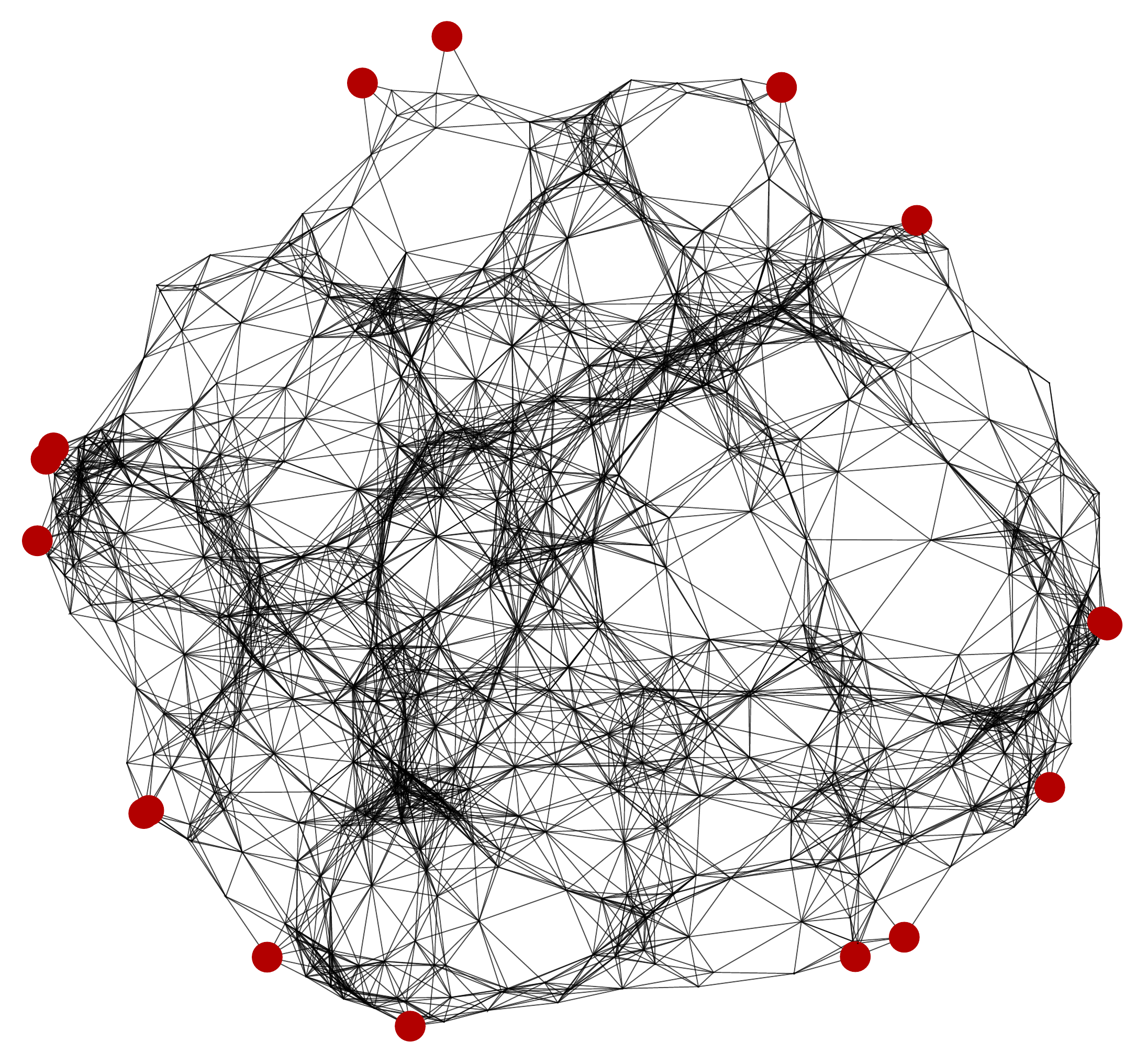}};
\end{tikzpicture}
\vspace{-10pt}
\caption{Two Gaussian Point Clouds with $\supp \mu$ in red.}
\end{figure}
\end{center}
\vspace{-10pt}

 In practice, the following two variations on the idea might be reasonable.
 \begin{enumerate}
\item  What is important for the properties of $\phi$ is not only the size of the support of $\mu$
but also the weights $\mu(w_i)$. We see from the explicit form of the embedding $\phi$ that if $\mu(w_i)$ is rather small, then the corresponding entries in the embedding
will vary very little. One could thus consider omitting and only focus on coordinates for which $\mu(w_i)$ is large. 
\item We know that $\phi$ sends $\supp \mu$ to a hyperplane and it might, in practice, make sense to move the embedding from $\mathbb{R}^m$ to
$\mathbb{R}^{m-1}$ either by PCA or by projecting onto the plane $x_1 + \dots + x_m = 0$.\\
\end{enumerate}

\textit{3. Swiss roll.} We will now apply both these ideas to an example of a manifold embedding: we generate a swiss roll in $\mathbb{R}^3$ using 10.000 points and building a graph connecting each points to its $k-$nearest neighbors (where $k=40$ but not tremendously important). Out of these 10.000 points, the measure ends up being supported on 50 points. Taking the vertices where $\mu$ is the largest, we realize that $1/3$ of the total probability mass is contained in only three vertices. We can take these three vertices to create an embedding into $\mathbb{R}^3$ which then, via PCA, we map to $\mathbb{R}^2$. Fig. 5 shows that the approximate shape is recovered by the embedding and that the embedding coordinates are smooth variables in the original space.
\vspace{-20pt}
\begin{center}
\begin{figure}[h!]
\begin{tikzpicture}
\node at (0,0) {\includegraphics[width=0.35\textwidth]{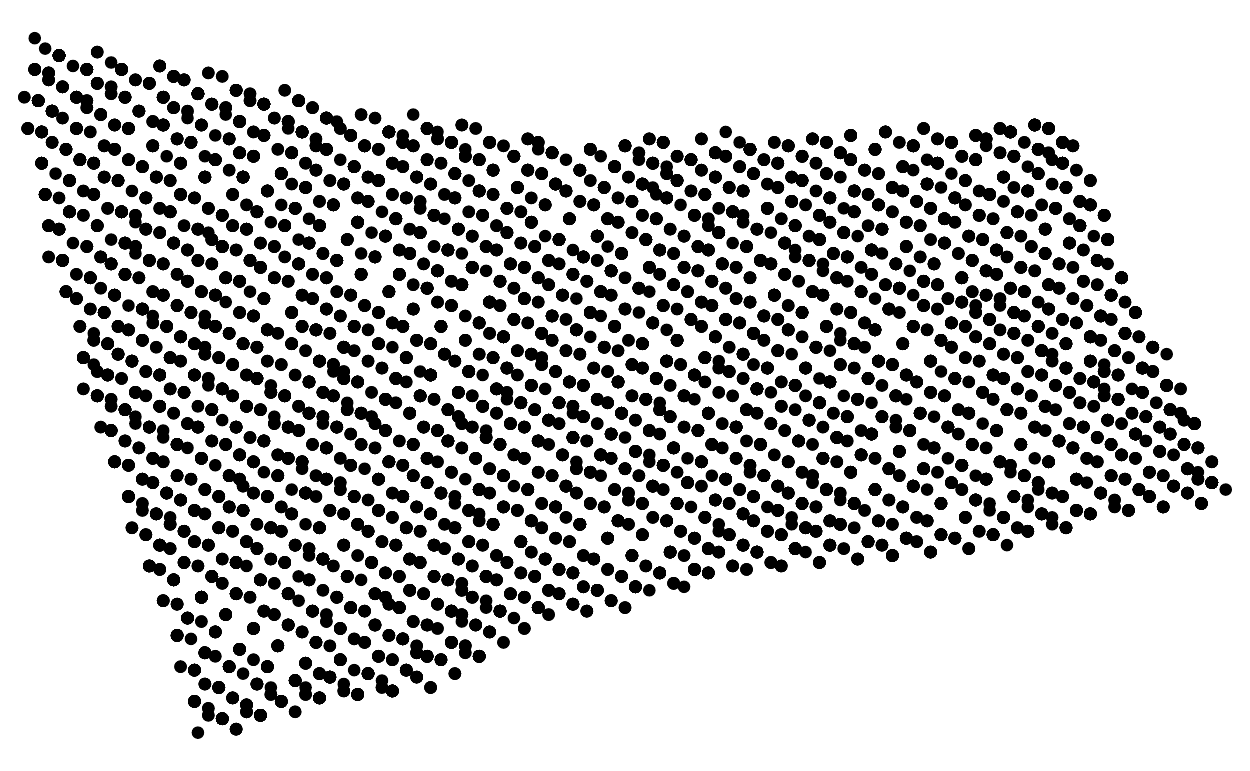}};
\node at (4,0) {\includegraphics[width=0.3\textwidth]{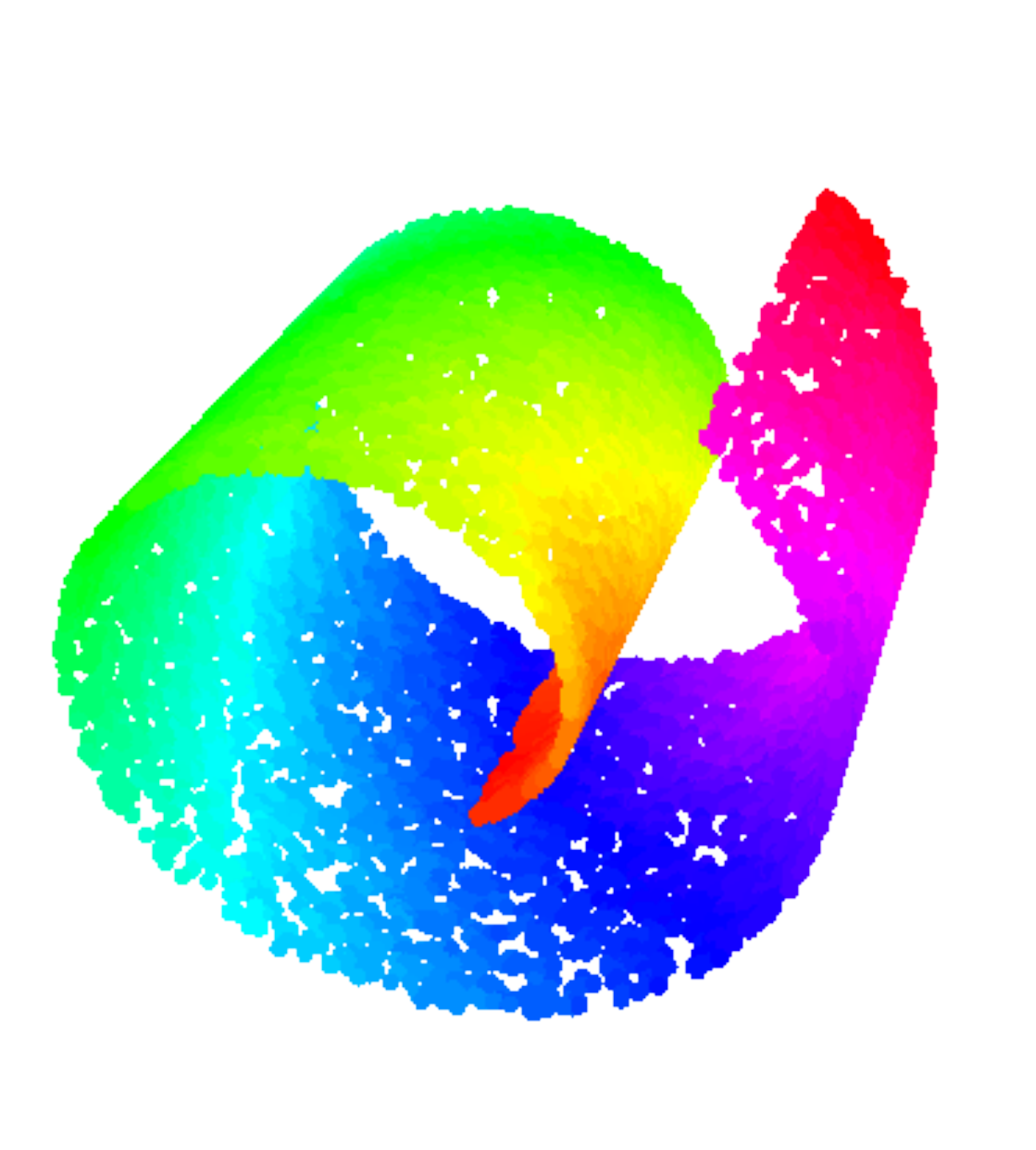}};
\node at (7.5,0) {\includegraphics[width=0.3\textwidth]{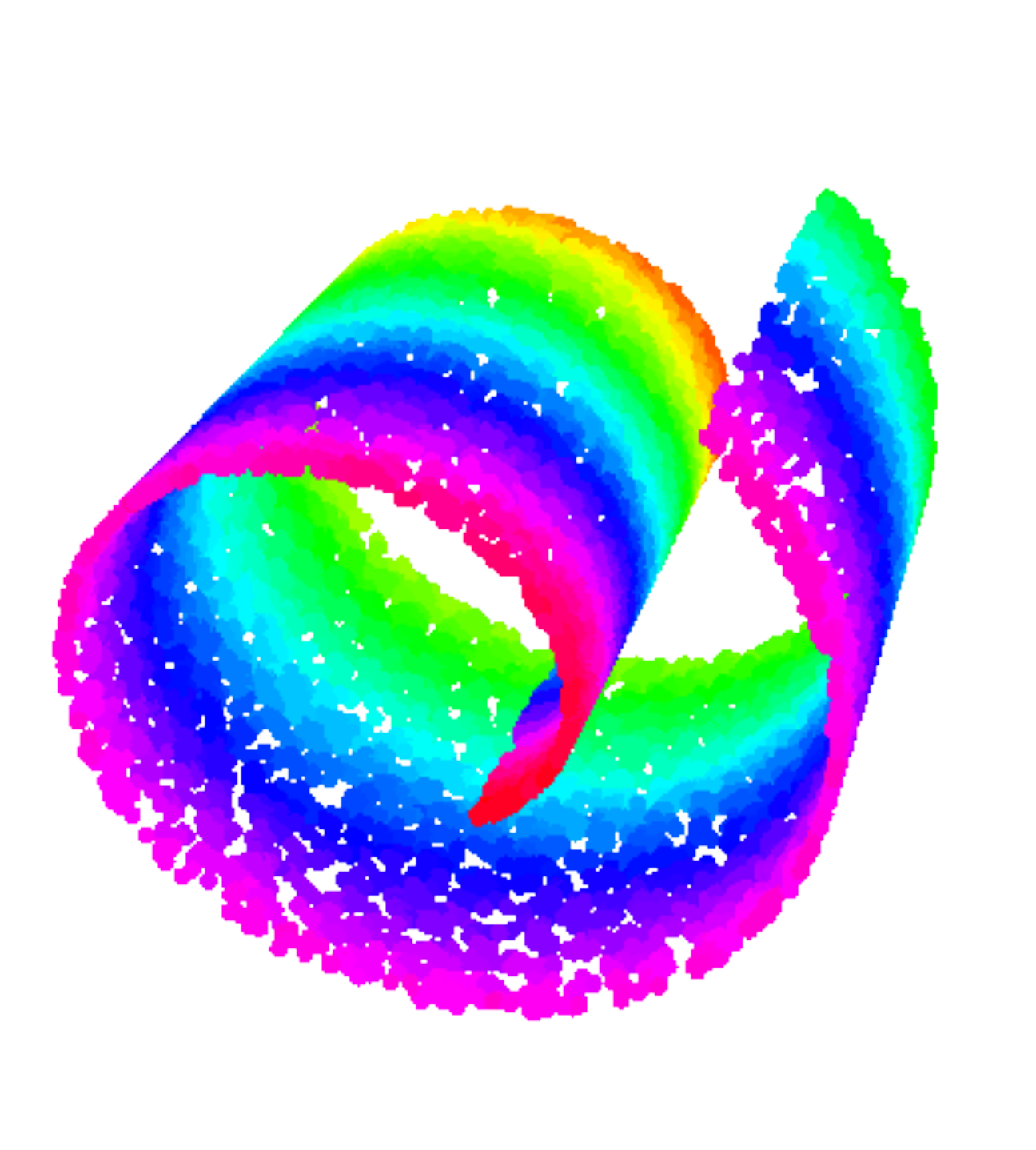}};
\end{tikzpicture}
\vspace{-10pt}
\caption{Left: embedding of a Swiss Roll into 2D using the $m=3$ followed by PCA. Original swiss role colored by $x-$coordinate (middle) and $y-$coordinate (right) of the embedding.}
\end{figure}
\end{center}
\vspace{-10pt}

\textit{4. Glued paths.} This interesting example was constructed by Noah Kravitz and is included with his kind permission. We take $m$ path graphs of length $2\ell+1$ and glue them together at the two endpoints.
This graph supports two very different balanced measures. One is concentrated at the two end-points:
the arising embedding is a line in $\mathbb{R}^2$ with diameter $\sim \diam(G)$ which collapses the
$m$ paths all into a singe path, it identifies the graph as being predominantly one-dimensional.

\begin{center}
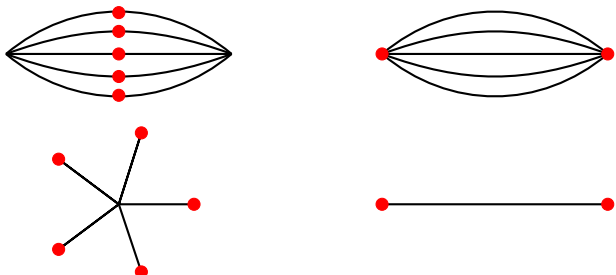
\begin{figure}[h!]
\begin{tikzpicture}
\draw [thick] (0,0) to[out=40, in=140] (3,0);
\draw [thick] (0,0) to[out=20, in=160] (3,0);
\draw [thick] (0,0) to[out=0, in=180] (3,0);
\draw [thick] (0,0) to[out=340, in=200] (3,0);
\draw [thick] (0,0) to[out=320, in=220] (3,0);
\filldraw[red] (1.5,0.55) circle (0.08cm);
\filldraw[red] (1.5,0.3) circle (0.08cm);
\filldraw[red] (1.5,0) circle (0.08cm);
\filldraw[red] (1.5,-0.55) circle (0.08cm);
\filldraw[red] (1.5,-0.3) circle (0.08cm);
\draw [thick] (5,0) to[out=40, in=140] (8,0);
\draw [thick] (5,0) to[out=20, in=160] (8,0);
\draw [thick] (5,0) to[out=0, in=180] (8,0);
\draw [thick] (5,0) to[out=340, in=200] (8,0);
\draw [thick] (5,0) to[out=320, in=220] (8,0);
\filldraw[red] (5,0) circle (0.08cm);
\filldraw[red] (8,0) circle (0.08cm);
\filldraw [thick] (1.5+1,-2+0) -- (1.5, -2) -- (1.5+0.3,-2+0.95) -- (1.5, -2) --(1.5-0.8,-2+0.6) -- (1.5, -2) --(1.5-0.8,-2-0.6)  -- (1.5, -2) ;
\filldraw [thick] (1.5+0.3,-2-0.9) -- (1.5, -2);
\filldraw[red] (1.5+1,-2+0) circle (0.08cm);
\filldraw[red] (1.5+0.3,-2+0.95) circle (0.08cm);
\filldraw[red] (1.5-0.8,-2+0.6) circle (0.08cm);
\filldraw[red] (1.5-0.8,-2-0.6) circle (0.08cm);
\filldraw[red] (1.5+0.3,-2-0.9) circle (0.08cm);
\filldraw [thick] (5, -2) -- (8,-2);
\filldraw[red] (5,-2) circle (0.08cm);
\filldraw[red] (8,-2) circle (0.08cm);
\end{tikzpicture}
\vspace{-10pt}
\caption{Taking $m$ path graphs of length $2\ell + 1$ all glued together at the endpoints: two balanced measures leading to two different embeddings emphasizing different aspects of the graph.}
\end{figure}
\end{center}

The other balanced measure enjoys a greater degree of stability and is equally concentrated on the $m$
midpoints of the $m$ paths. We see that, for this example, Theorem 2 is sharp up to constants: for all $v \in \supp \mu$
$$ \frac{1}{\# \supp \mu} \sum_{ w \in \supp \mu} \| \phi(v) - \phi(w) \|_{\ell^{1}(\mathbb{R}^m)} \sim  \frac{\diam(G)}{m}.$$
In particular, the graph is being folded into a rather small region of $\ell^1(\mathbb{R}^m)$ whose diameter shrinks as
$m$ increases. Nonetheless, the embedding itself is certainly a very good representation of the structure of the graph.

\section{Comments, Examples and Related Results}

\subsection{Balanced measures as critical points.}
We recall that for any given graph $G=(V,E)$, we call a probability measure $\mu$ on the set of vertices balanced if
$$ \mu(u) > 0 \implies \sum_{v \in V} d(u,v) \mu(v) = \max_{w \in W} \sum_{v \in V} d(w, v) \mu(v).$$
There is a simple variational characterization of balanced measures. The distance matrix $D \in \mathbb{R}^{n \times n}$ defined via $D_{ij} = d(v_i, v_j)$, has integer
entries and is symmetric. Moreover, we identify the set $\mathcal{P}(V)$ of probability measure on $V$ with 
$$ \Delta = \left\{(x_1, \dots, x_n) \in \mathbb{R}_{\geq 0}^{n}: x_1 + \dots + x_n = 1\right\}.$$
\begin{proposition}
Each critical point of the functional $J:\Delta \rightarrow \mathbb{R}_{\geq 0}$ given by
$$ J(\mu) = \left\langle \mu, D \mu \right\rangle$$
is a balanced measure. A balanced measure $\mu$ has the property that the directional derivative in all admissible directions is non-positive.
\end{proposition}
If $\mu$ is on the boundary, we say it is a critical point if all directional derivatives in admissible directions vanish: for any
signed measure $\nu$ such that $\mu + \varepsilon \nu$ is a probability measure for all $0 \leq \varepsilon \leq \varepsilon_0$
for some $\varepsilon_0 > 0$ (depending on $\mu$), we require 
$$J(\mu+ \varepsilon \nu) = J(\mu) + o(\varepsilon) \qquad \mbox{as} ~\varepsilon \rightarrow 0^+.$$
Since $\Delta$ is a compact set and $J$ is continuous, there is at least one maximum and thus there is
always at least one balanced measure. We note that the global minima are given by the Dirac measures in a 
single vertex: these are \textit{not} critical points (possible because they are assumed at the boundary of $\Delta$).

\subsection{Balancing measures and boundary}  One could wonder whether there is any way of deciding a priori where the measure $\mu$ can be supported. 
 We recall a definition of \textit{boundary} of a graph given in \cite{stein1}: for any connected graph $G=(V,E)$, we define
the boundary $\partial G \subseteq V$ to be the set of all vertices $u \in V$ for which there exists another vertex $v$ such that the \textit{average} neighbor of $u$ is
closer to $v$ than $d(u,v)$. Formally, we define the boundary as
$$\partial G = \left\{u \in V \big|  ~\exists v \in V:  ~  \frac{1}{\deg(u)} \sum_{(u, w) \in E} d(w,v) < d(u,v)  \right\}.$$

\begin{center}
\begin{figure}[h!]
\begin{tikzpicture}
\node at (0,0) {\includegraphics[width=0.27\textwidth]{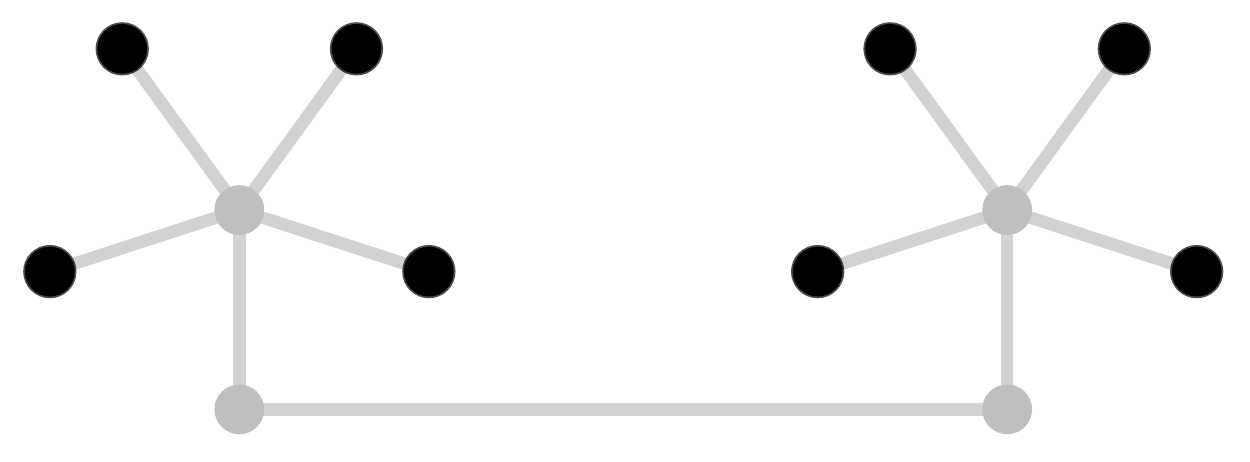}};
\node at (3.25,0) {\includegraphics[width=0.1\textwidth]{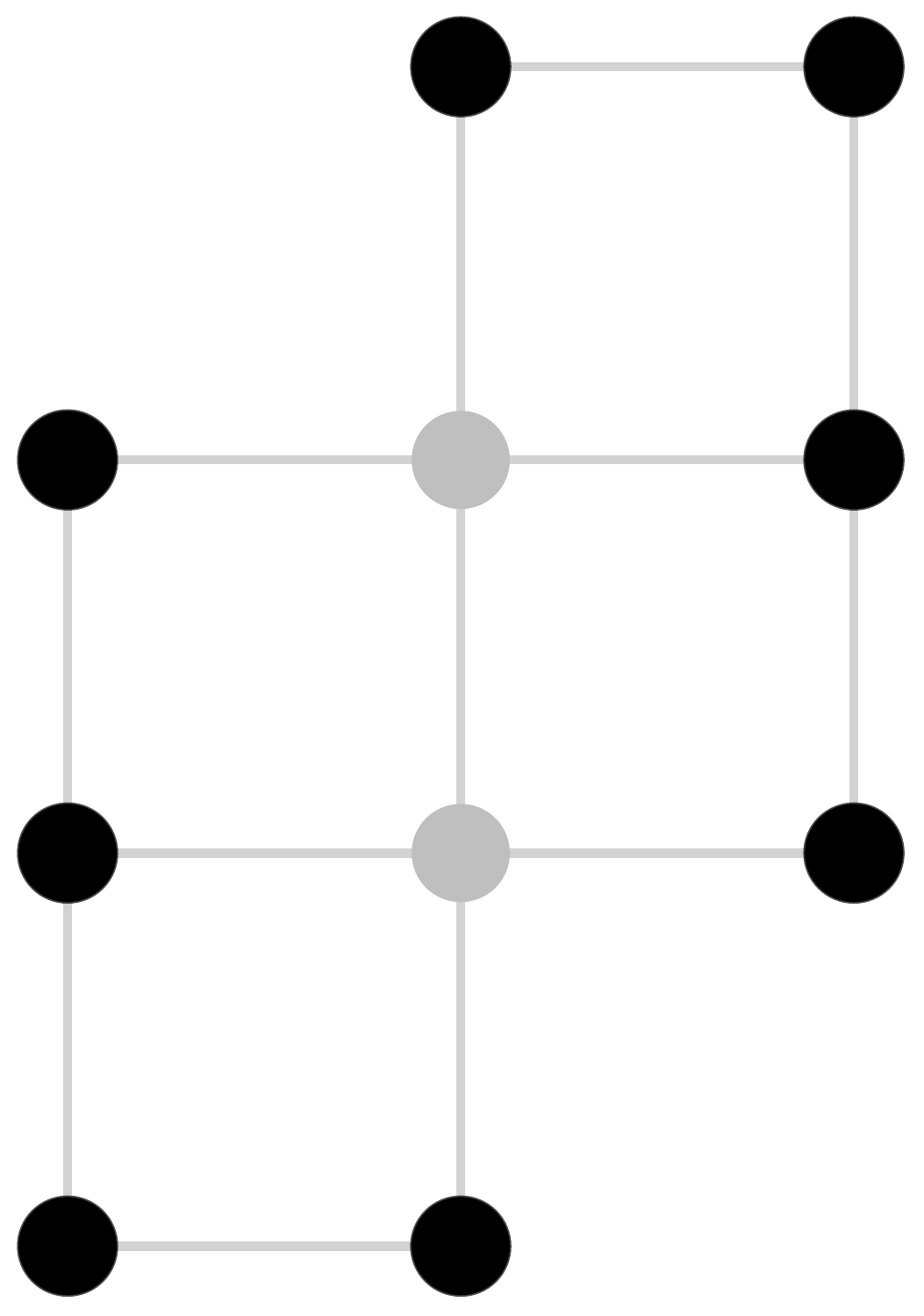}};
\node at (6.5,0) {\includegraphics[width=0.18\textwidth]{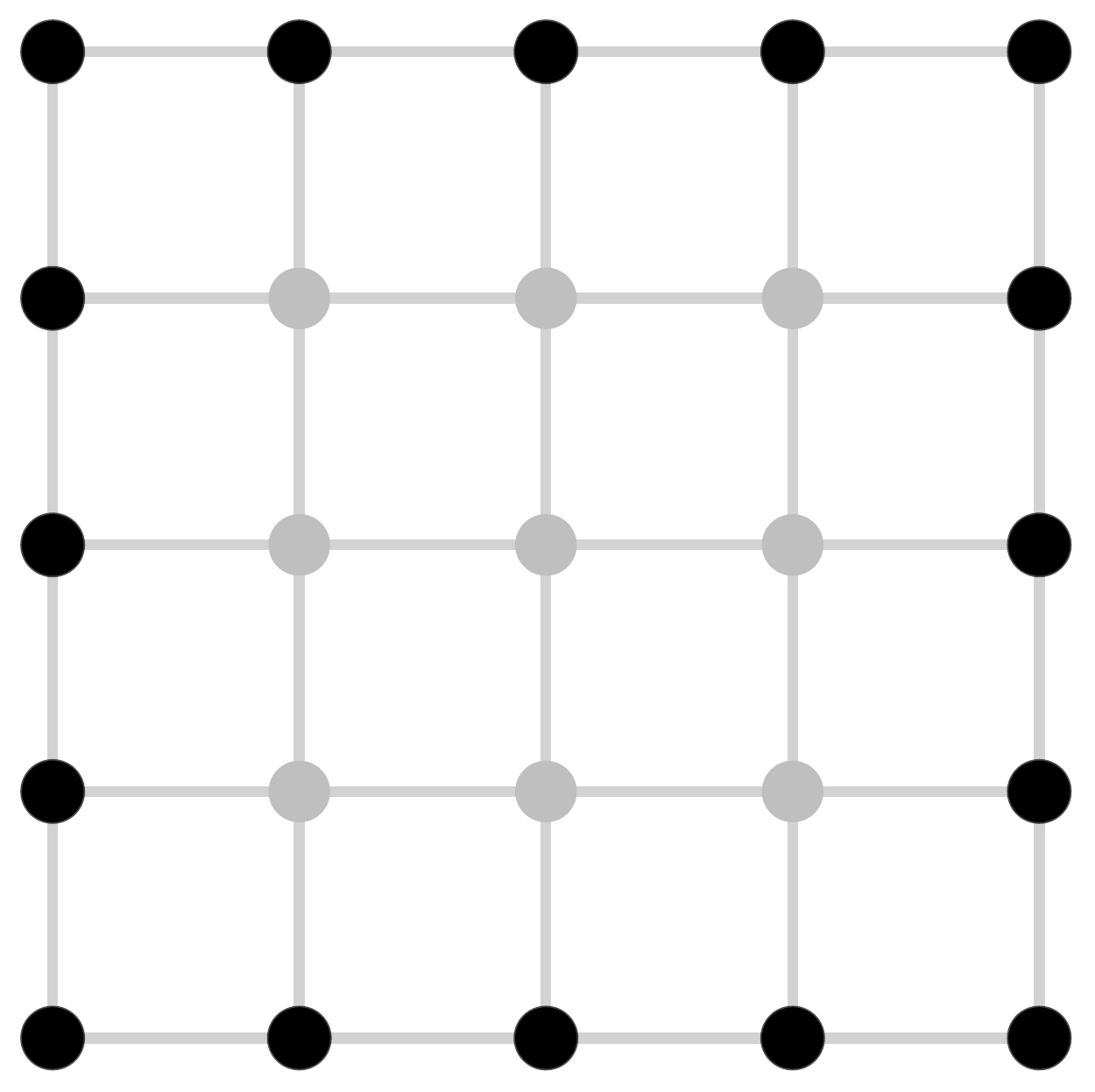}};
\end{tikzpicture}
\vspace{-10pt}
\caption{Three graphs with their boundary vertices highlighted.}
\end{figure}
\end{center}

This notion of boundary is motivated by the fact that it satisfies a type of isoperimetric inequality stating that large graphs have to have 
a large boundary: if the maximal degree of $G$ is given by $\Delta$, then
$$ \# \partial G \geq \frac{1}{2\Delta} \frac{\#V}{\diam(G)}.$$
We refer to \cite{stein1} for more details. As it turns out, we will be able to show that the long-time behavior of the greedy procedure tends to happen in the boundary.

\begin{proposition}
Let $x_1, \dots, x_k \in V$ be a list of vertices. Then
$ f(v) = \sum_{i=1}^{k} d(v,x_i)$ assumes its maximum in $\partial G.$
If $f$ assumes its maximum in $V \setminus \partial G$, $f$ is constant. In particular, 
there exists a balanced measure supported in $\partial G$.
\end{proposition}
One might assume that this means that `typically' balanced measures can only be supported in
the boundary -- it would be interesting to understand this better.

\subsection{Remarks.} 
 In many of the examples shown throughout the paper, we see that $\# \supp(\mu) = m \ll n$. This need not always be the case: examples are given by the dodecahedral graph
or the Desargue graphs for which the uniform measure is balanced. 
Both examples have distance matrices $(d(v_i, v_j))_{i,j=1}^n$ with one positive eigenvalue whose corresponding
eigenvector is constant which, in light of Proposition 1, is perhaps not a coincidence. In practice, it seems very difficult for graphs to have $\supp \mu$ contain a large number of vertices and it would be interesting to have a more quantitative understanding of this. Are there upper bounds on $\# \supp \mu$ depending on some graph parameters?\\

 Graphs with symmetries will naturally present with balanced measures inheriting these symmetries. However, the support of a balanced measure is \textit{not} a graph invariant: graphs can have balanced measures with very different cardinality (see Fig. 7). It would be interesting to understand how many different balanced measures a graph can support.
 A graph having a balanced measure supported in two vertices indicates that the graph is elongated and behaves, mostly, like a one-dimensional interval -- there are many graphs like this. Conversely, a graph supporting a balanced measure in $m \geq n/100$ vertices seems to require a tremendous amount of symmetries and these are probably  rare. It might also be interesting to study the properties of balanced measures on Erd\H{o}s-Renyi random graphs. \\
\vspace{-20pt}
\begin{center}
\begin{figure}[h!]
\begin{tikzpicture}
\node at (0,0) {\includegraphics[width=0.32\textwidth]{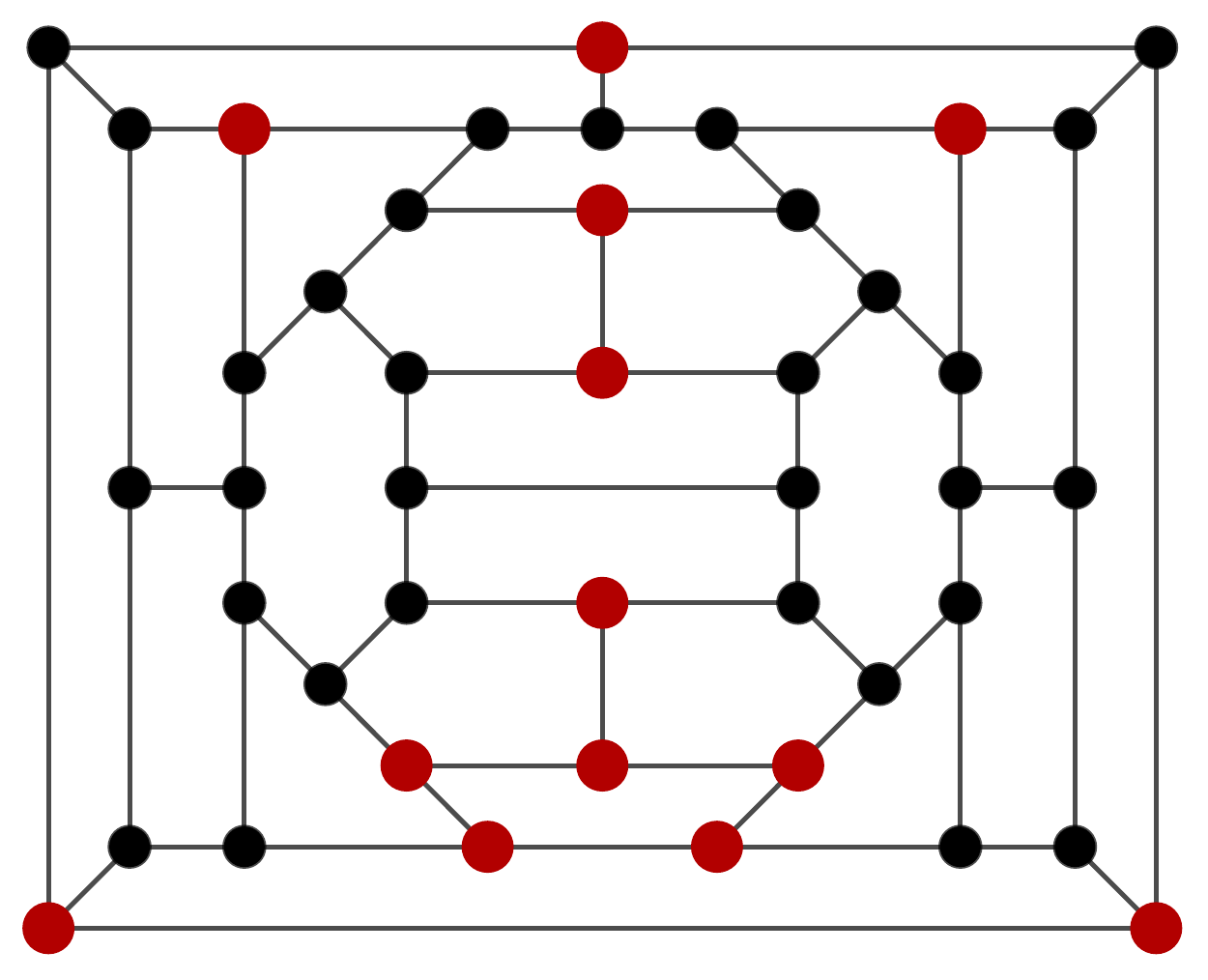}};
\node at (5,0) {\includegraphics[width=0.32\textwidth]{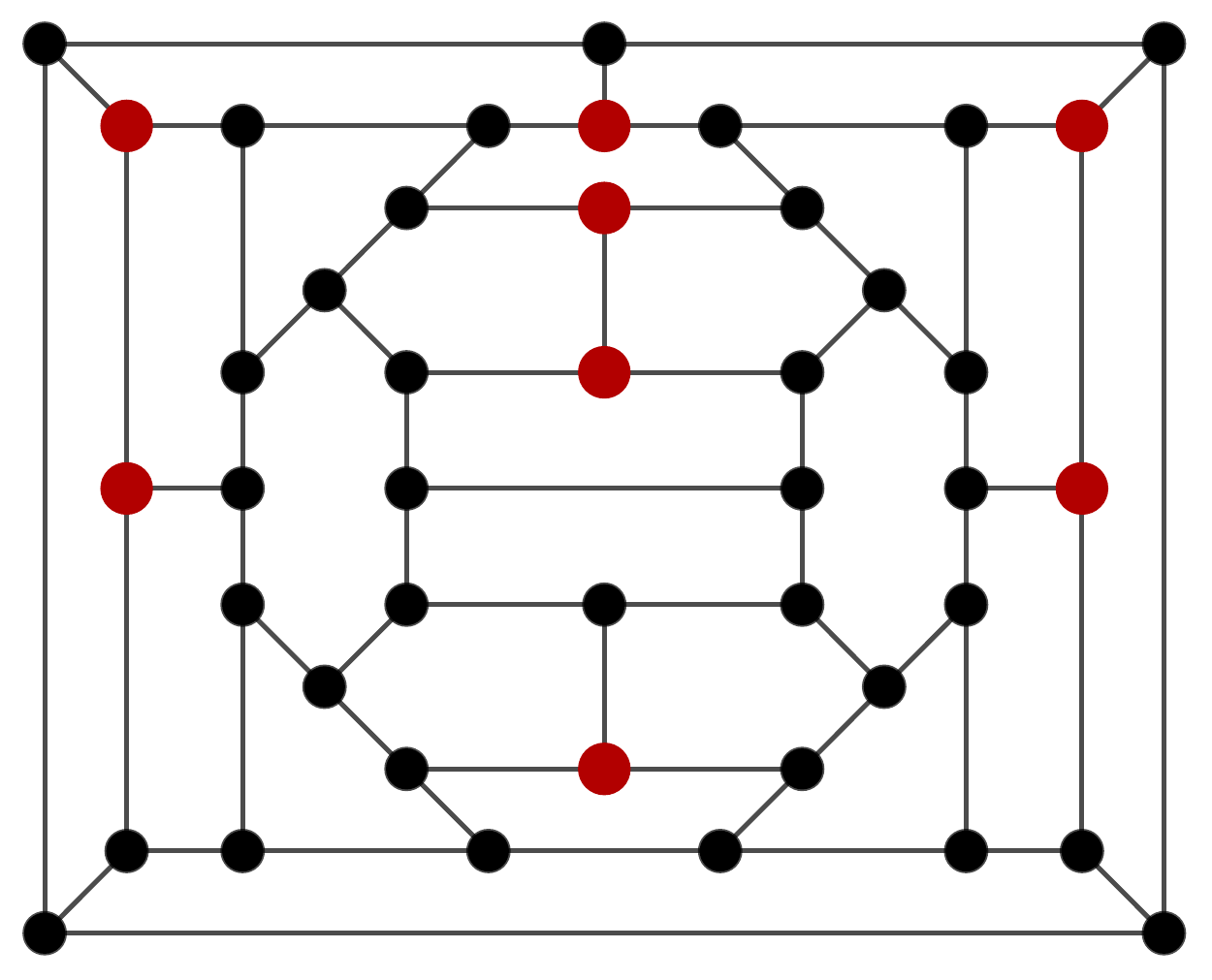}};
\end{tikzpicture}
\vspace{-10pt}
\caption{Support of two balanced measures on Grinberg 44.}
\end{figure}
\end{center}
\vspace{-10pt}

We conclude with an interesting phenomenon that appears to be rare. Recall that the definition of a balanced measures is given by
$$ \mu(u) > 0 \implies \sum_{v \in V} d(u,v) \mu(v) = \max_{w \in W} \sum_{v \in V} d(w, v) \mu(v).$$
One could wonder whether any type of inverse result of this flavor is true: are the maxima \textit{only} assumed
in vertices where the measure is supported? The Frucht graph (see Fig. 1) is a counterexample: there is a balanced measure on 4 vertices (with weights $0.1, 0.2, 0.3, 0.4$) such that $\sum d(v,w) \mu(w)$
assumes its maximum in 6 vertices. However, this seems to be exceedingly rare.

%\begin{center}
%\begin{figure}[h!]
%\begin{tikzpicture}
%\node at (0,0) {\includegraphics[width=0.6\textwidth]{ex5}};
%\end{tikzpicture}
%\vspace{-15pt}
%\caption{Thomassen 94.}
%\end{figure}
%\end{center}
%\vspace{-20pt}

\subsection{Related results.} \textit{Potential Theory.} Let $(X,d)$ be a compact metric space and let $\mu$ be a probability measure on the
space. The study of the energy integral
$$ J(\mu) = \int_{X \times X}  d(x,y) d\mu(x)  d\mu(y)$$
is classical and dates back to a 1956 paper of Bj\"orck \cite{bj}. Bj\"orck proved that if $X$ is a compact subset of $\mathbb{R}^n$,
then the maximizing measure is unique and supported on the boundary $\partial X$. This has inspired a lot of work
in Euclidean space: we refer to Alexander \cite{alex, alex2}, Alexander \& Stolarsky \cite{alex3}, 
Carando, Galicer \& Pinasco \cite{carando}, Hinrichs, Nickolas \& Wolf \cite{hinrichs} and Wolf \cite{wolf2}.
A continuous version of the greedy procedure was proposed by Larcher, Schmid \& Wolf \cite{larcher} to approximate the
(unique) maximizing measure. Outside the Euclidean setting there has been substantial work in quasihypermetric spaces, we refer to Nickolas \& Wolf \cite{nick0, nick1, nick2, nick3, wolf}. The Euclidean space is an example of a quasihypermetric space; in our setting, a graph is quasihypermetric if the distance matrix $D$ is negative semi-definite when restricted to vectors whose sum of entries is 0.
Partial motivation for many of these questions is a fascinating 1964 result of Gross \cite{gross}: in any compact, connected metric space $(X,d)$, there exists a real number $r(X,d) \in \mathbb{R}$ such that for any set of points $\left\{x_1, \dots, x_n\right\} \subset X$, there exists a point $x \in X$ with
$$ \frac{1}{n} \sum_{i=1}^{n} d(x, x_i) = r(X,d).$$
A wonderful introduction to the subject is given by Cleary, Morris \& Yost \cite{cleary}.
No such results can be true on graphs, however, see Thomassen \cite{thomassen} for a substitute result (see also \cite{stein2} for a connection to the
von Neumann Minimax Theorem). The result of Gross can be extended to probability measures which has the following implication: if there exists a probability measure $\mu$ on $X$ so that
$$ \int_{X} d(x,y) d\mu(y) \qquad \mbox{is independent of}~x,$$
then $r(X,d)$ has to equal the value of the integral evaluated at an arbitrary point $x$. This allows for an explicit computation of $r(X,d)$ which, in general,
is very difficult. Simultaneously, if we consider the problem of maximizing
$$ \mu \rightarrow \int_{X \times X}  d(x,y) d\mu(x)  d\mu(y),$$
then an Euler-Lagrange ansatz shows that $\int_{X} d(x,y) \mu(y)$ is constant on $\supp(\mu)$. If $\supp(\mu) =X$, the maximizing measure can be used to determine $r(X,d)$. Our paper is somewhat opposite 
to the main direction in the literature since (besides the graph setting which is a metric space but \textit{not} a connected metric space) we are interested in cases where $\supp(\mu)$ is much smaller than the space $X$.\\
\textit{Beacons.}
Another related area is in the theory of networks. A \textit{beacon-based} approach is to select, often randomly, a number of vertices (the `beacons'). One is then given the corresponding rows in the Graph Distance Matrix, the distance between the beacons and all other vertices and tries to reconstruct all distances. An example of such a result is given by Ng \& Zhang \cite{ng} who showed that a beacon-based approach could embed all but a small fraction of Internet distances; see also Kleinberg, Slivkins \& Wexler \cite{kleinberg}. These results seem are somewhat different in flavor insofar as the beacons can be freely chosen whereas in our setting they arise from balanced measures which are determined by the graph.\\
\textit{Greedy sequences.} Another motivation of the present paper are the recent introduction of greedy sequences for the purpose of generating uniformly distributed sequences with good regularity properties \cite{brown,kritzinger,steind,steind2,steind3}. This is in stark contrast to the present paper where the limiting measure $\mu$ tends to be highly irregular, however, this is due to the particular type of kernel being used; kernels with the purpose of generating a more uniform distribution on graphs have been studied by Brown \cite{brown2}, Cloninger \& Mhaskar \cite{alex} and Linderman \& Steinerberger \cite{george}.

\section{Proofs}

\subsection{Proof of Theorem 1}
The proof decouples into the following steps.
\begin{enumerate}
\item 
$$ E_m = \frac{1}{m(m-1)}\sum_{i,j=1}^{m} d(x_i, x_j) \qquad \mbox{converges as} ~m \rightarrow \infty$$
\item and
$$ \frac{\diam(G)}{2} \leq \lim_{m \rightarrow \infty} E_m \leq \diam(G).$$
\item A weighted average of $T_m(v) = \sum_{v \in V} d(v,w) \mu_m(w)$ is close to $E_m$.
\item A stability estimate: if $T_m(v)$ is large, then so is $T_{m + \ell}(v)$ for small $\ell \in \mathbb{N}$.
\item From this we deduce an upper bound on $\|T_m\|_{\ell^{\infty}}$.
\item Steps (3) and (5) imply that $T_m(v)$ is nearly constant on $\supp \mu$.
\end{enumerate}
\begin{proof}  \textbf{Step 1.} We start by considering the behavior of the rescaled energy
$$ E_m = \frac{1}{m(m-1)}\sum_{i,j=1}^{m} d(x_i, x_j).$$
The first step of the argument consists in showing that $E_m$ converges.
This sequence is trivially bounded since
$$ E_m = \frac{1}{m(m-1)}\sum_{i,j=1}^{m} d(x_i, x_j) \leq \frac{1}{m(m-1)}\sum_{i,j=1 \atop i \neq j}^{m} \diam(G) = \diam(G).$$
We observe the algebraic identity
$$ (m+1)m \cdot E_{m+1} = m (m-1)  \cdot E_m + 2 \sum_{i=1}^{m} d(x_{m+1}, x_i).$$
Recalling the definition of the greedy sequence
$$ \sum_{i=1}^{m} d(x_{m+1}, x_i) = \max_{v \in V} \sum_{i=1}^{m} d(v, x_i).$$
We can estimate this from below by replacing the maximum with any weighted average: for any probability measure $\nu$ on $V$, we have
$$ \max_{v \in V} \sum_{i=1}^{m} d(v, x_i) \geq \sum_{v \in V} \sum_{i=1}^{m} d(v, x_i) \nu(v).$$
Our choice will be the empirical measure of the first $m$ points
$$ \nu  = \frac{1}{m} \sum_{j=1}^{m} \delta_{x_j}.$$
For this particular choice of probability measure, the lower bound simplifies 
\begin{align*}
  \max_{v \in V} \sum_{i=1}^{m} d(v, x_i)  &\geq \sum_{v \in V} \sum_{i=1}^{m} d(v, x_i) \nu(v) \\
  &= \frac{1}{m} \sum_{i,j=1}^{m} d(x_i, x_j) = (m-1) \cdot E_m.
  \end{align*}
Plugging in shows that this implies
 $$ (m+1)m E_{m+1} \geq (m+2)(m-1)  E_m$$
and thus
$$ E_{m+1} \geq \frac{(m+2)(m-1)}{m (m+1)} \cdot E_m.$$
An expansion shows that, as $m \rightarrow \infty$,
$$ \frac{(m+2)(m-1)}{m (m+1)} \geq 1 - \frac{2}{m^2}.$$
This implies, in particular, that for all $m$ sufficiently large and all $l \in \mathbb{N}$,
$$ E_{m+\ell} \geq E_m \cdot \prod_{k=m}^{m+l-1} \left(1 - \frac{2}{k^2}\right).$$
This product converges as $\ell \rightarrow \infty$ to a number that is arbitrarily close to 1
as $m \rightarrow \infty$. We will now make this quantitative.
Using standard Taylor series estimates, we have for $N$ sufficiently large 
\begin{align*}
\prod_{m=N}^{\infty} \left(1 - \frac{2}{m^2}\right) &= \exp \left[\log  \left(\prod_{m=N}^{\infty} \left(1 - \frac{2}{m^2}\right) \right) \right] \\
&= \exp \left[ \sum_{m=N}^{\infty} \log \left(1 - \frac{2}{m^2}\right) \right]  \\
&\geq \exp  \left( \sum_{m=N}^{\infty} - \frac{4}{m^2} \right) \geq \exp \left(-\frac{8}{N}\right) \geq 1 - \frac{16}{N}.
\end{align*}
This implies that, for $m$ sufficiently large,
$$ \inf_{k \geq m} E_k \geq \left(1 - \frac{16}{m}\right)E_m.$$
Thus, given $E_m$, subsequent values cannot be much smaller. Let 
$$ \alpha =  \limsup_{m \rightarrow \infty} E_m.$$
Then, for any $\varepsilon > 0$ there exists a subsequence $(E_{m_{\ell}})_{\ell=1}^{\infty}$ such that for all elements of the sequence $E_{m_{\ell}} \geq \alpha - \varepsilon/2$. Picking $\ell$ sufficiently large so that $16/m_{\ell} \leq \varepsilon/(2\alpha)$, we deduce that
$$ \liminf_{m \rightarrow \infty} E_m \geq \alpha - \varepsilon$$
and since $\varepsilon$ was arbitrary, we deduce the convergence of $(E_m)_{m=1}^{\infty}$. \\

\textbf{Step 2.} We denote the limit by
$$ \alpha = \lim_{m \rightarrow \infty} E_m$$
and will now establish $\diam(G)/2 \leq \alpha \leq \diam(G)$. The upper bound is obvious, it remains to prove the lower bound.
At this point we invoke the von Neumann Minimax Theorem in the special case of symmetric matrices (see \cite{stein2, thomassen}) and apply it to
the case where the matrix is given by the distance matrix $D_{ij} = d(v_i, v_j) \in \mathbb{R}^{n \times n}$. It implies that there exists a constant $\beta$ depending only on the matrix
such that for all probability vectors $\mu$ 
$$ \min_{1 \leq i \leq n} (D\mu)_i \leq \beta \leq  \max_{1 \leq i \leq n} (D\mu)_i.$$
We now choose a specific measure. By picking two vertices $a,b \in V$ at maximal distance $d(a,b) = \diam(G)$ and then choosing
$$ \nu = \frac{1}{2} (\delta_a + \delta_b),$$
we can use the triangle inequality to observe that for any other vertex $v$
$$ \frac{\diam(G)}{2} = \frac{d(a,b)}{2} \leq \frac{d(a,v) + d(v,b)}{2} = (D\nu)(v).$$
This implies that $\beta \geq \diam(G)/2$. Thus
\begin{align*}
 (m+1)m \cdot E_{m+1} &= m (m-1)  \cdot E_m + 2 \sum_{i=1}^{m} d(x_{m+1}, x_i) \\
 &\geq m (m-1)  \cdot E_m + m \cdot \diam(G).
 \end{align*}
Abbreviating $F_m = m (m-1) E_m$, we have $F_{m+1} \geq F_m + m \cdot \diam(G)$ and thus
$$ F_{m} \geq \diam(G) \sum_{k=1}^{m-1} k = (1+o(1)) \cdot \diam(G) \cdot \frac{m^2}{2}$$
which implies the desired statement.\\

\textbf{Step 3.} Recalling $\alpha = \lim_{m \rightarrow \infty} E_m$, we can use the inequality 
$$ \inf_{k \geq m} E_k \geq \left(1 - \frac{16}{m}\right)E_m$$
in combination with $\alpha \geq  \inf_{k \geq m} E_k$ to conclude, for $m$ sufficiently large,
$$ E_m \leq \left(1 + \frac{20}{m}\right) \alpha.$$
These two facts combined suggest that we would except $E_m$ to approach its limit from below (though this will not be used in subsequent arguments).
Let now $\varepsilon > 0$ be arbitrary and let $m$ be so large that
$$ \forall~k \geq m: \quad \quad E_k \geq \alpha - \varepsilon.$$
We introduce the (transport cost) function $T_m: V \rightarrow \mathbb{R}$ via
$$ T_m(v) = \sum_{w \in V} d(v,w) \mu_m(w).$$
We observe that the $\mu_m-$weighted average of $f_m$ is explicit since
\begin{align*}
 \sum_{v \in V}  \mu_m(v) T_m(v) &=  \sum_{v,w \in V}  \mu_m(v) d(v,w) \mu_m(w) \\
 &= \frac{1}{m^2}   \sum_{i,j=1}^m d(x_i,x_j) = \frac{m-1}{m} E_m
\end{align*}
which, for $m$ very large, is very close to $E_m$ which in turn is close to $\alpha - \varepsilon$. We will now show that
the maximal value of $T_m(v)$ is also close to $\alpha$.\\

\textbf{Step 4.} The basic ingredient is a continuity estimate: if $\|T_m\|_{\ell^{\infty}}$ was much larger than $\alpha$, then we would expect
this to also be true for $T_{m+1}$ since $\mu_{m+1}$ is rather similar to $\mu_m$. Indeed, we have
$$ \mu_{m+1} = \frac{1}{m+1} \sum_{k=1}^{m+1} \delta_{x_k} = \frac{m}{m+1} \mu_m + \frac{1}{m+1} \delta_{x_{m+1}}$$
and therefore, since $T_m \geq 0$ and $T_m(v) \leq \diam(G)$,
\begin{align*}
 \|T_{m+1} - T_m\|_{\ell^{\infty}} = \left\| \frac{m}{m+1} T_m + \frac{d(x_{k+1},v)}{m+1}  - T_m\right\|_{\ell^{\infty}} \leq \frac{\diam(G)}{m+1} .
 \end{align*}
Suppose now there exists $v \in V$ such that
$$ T_m(v) \geq \alpha + \delta$$
for some $\delta > 0$. Then the value of $T_{m+1}, T_{m+2}, \dots$ at the same vertex is still large at least for the next few iterations of the process and we have for all $\ell \in \mathbb{N}$
$$  \|T_{m+\ell}\|_{\ell_{\infty}} \geq T_{m+\ell}(v) \geq \alpha + \delta - \frac{\ell}{m+1} \diam(G)$$
and thus, for 
$$1 \leq \ell \leq \frac{m \delta}{4 \diam(G)},$$ 
we have that 
$$T_{m+\ell}(v) \geq \alpha + \delta/2.$$

\textbf{Step 5.} Let $\varepsilon > 0$ fixed and all $m$ be so large 
$$  \forall~k \geq m \qquad \quad \alpha - \varepsilon \leq E_k \leq \left(1 + \frac{20}{k} \right)\alpha.$$
This will now be shown to imply that
$$ \| T_m\|_{\ell^{\infty}} \leq \alpha + 3 \sqrt{\diam(G)} \sqrt{\varepsilon}$$
Let us suppose that, for some $\delta>0$, we have $T_m(v) \geq \alpha + \delta$. We use Step 4 with $\ell = m \delta/ (4 \diam(G))$ to argue that since
in each step we add at least $\alpha + \delta/2$, 
\begin{align*}
 ( m+ \ell) (m+\ell -1) E_{m+\ell} &\geq m(m-1) E_m + 2 \ell m\left(\alpha + \frac{\delta}{2} \right) \\
 &\geq m(m-1) (\alpha - \varepsilon) + 2\ell m\left(\alpha + \frac{\delta}{2} \right).
 \end{align*}
 We will now bound the upper term from above and the lower term from below. For the bound from above, we first
 recall that, for $m$ sufficiently large,
$$ E_{m+\ell}  \leq \left(1 + \frac{20}{m}\right) \alpha.$$
We bound from above, using $\delta \leq \diam(G)$,
\begin{align*}
  ( m+ \ell) (m+\ell -1) E_{m+\ell}  &\leq   \left(1 + \frac{20}{m}\right)\alpha (m+\ell)^2 \\
  &=  \left(1 + \frac{20}{m}\right)\alpha m^2 \left(1 + \frac{\delta}{4 \diam(G)}\right)^2 \\
  &\leq \alpha m^2 \left(1 + \frac{\delta}{4 \diam(G)}\right)^2 +  \frac{20}{m} \alpha m^2 \left(1 + \frac{\delta}{4 \diam(G)}\right)^2\\
  &\leq 40 \alpha m + \alpha m^2 \left(1 + \frac{\delta}{4 \diam(G)}\right)^2.
  \end{align*}
Using $\alpha \leq \diam(G)$, this leads to
  \begin{align*}
  ( m+ \ell) (m+\ell -1) E_{m+\ell}  \leq  40 \alpha m + \alpha m^2 + \frac{ \alpha \delta m^2}{2\diam (G)} + \frac{m^2 \delta^2}{16 \diam(G)} .
  \end{align*}
 For the lower bound, we argue that by plugging in the definition of $\ell$ we have
  \begin{align*}
   m(m-1) (\alpha - \varepsilon) + 2\ell m\left(\alpha + \frac{\delta}{2} \right) &\geq (m-1)^2 (\alpha - \varepsilon) + \frac{m^2 \delta \alpha}{2 \diam(G)} + \frac{m^2 \delta^2}{4\diam(G)}.
  \end{align*}
  We see that this leads to a contradiction as soon as
  $$  (m-1)^2 \varepsilon + 40 \alpha m \leq \frac{m^2 \delta^2}{8 \diam(G)}$$
  and thus, for $m$ sufficiently large,
  $$ \delta \leq 3 \sqrt{\diam(G)} \sqrt{\varepsilon}.$$
  This shows that 
  $$ \limsup_{m \rightarrow \infty} \| T_m\|_{\ell^{\infty}} \leq \alpha.$$
  
  \textbf{Step 6.} We will now refine this last statement and show that $\limsup$ can be replaced by a $\lim$ on $\supp \mu_m$. Combining the estimates from the previous sections 
    \begin{align*}
  \frac{m-1}{m} E_m &=  \frac{1}{m^2}   \sum_{i,j=1}^m d(x_i,x_j) \\
  &=  \sum_{v,w \in V} \mu_m(v) d(v,w) \mu_m(w) \\
  &=  \sum_{v \in V} T_m(v) \mu_m(v) .
\end{align*}
Recalling that $\mu_m$ is a probability measure, that
$$ \lim_{m \rightarrow \infty} E_m = \alpha$$
as well as
$$ \lim_{m \rightarrow \infty} \|T_m\|_{\ell^{\infty}} = \alpha,$$
we conclude that $T_m$ has to be close to maximal in each vertex carrying a non-vanishing portion of the probability mass in the sense of (3) of the statement.
\end{proof}

\subsection{Proof of Theorem 2}
\begin{proof}
Let $G=(V,E)$ be arbitrary and let $\mu$ be a balanced measure supported on $m$ vertices $w_1, \dots, w_m \in V$. We define the map $\phi: V \rightarrow \mathbb{R}_{\geq 0}^m$ via
$$ \phi(v) = \left( \mu(w_1) d(v, w_1), \dots, \mu(w_m) d(v, w_m) \right).$$
We first note that $\phi$ maps $V$ into the positive orthant $\mathbb{R}^m_{\geq 0}$. Moreover, we observe that since
$$ T(v) = \sum_{w \in V} d(w,v) \mu(v) \qquad \mbox{assumes its maximum in the support},$$
we have, for each $1 \leq i \leq m$,
\begin{align*}
\| \phi(v)\|_{\ell^1} = \sum_{j=1}^{m}  d(v, w_j) \mu(w_j) \leq \sum_{j=1}^{m} d(w_i, w_j) \mu(w_j).
\end{align*}
Summing over $i$, we arrive at
$$ \| \phi(v)\|_{\ell^1}  \leq  \sum_{i=1}^{m} \mu(w_i) \sum_{j =1}^{m} d(w_i, w_j)  \mu(w_j) = \alpha.$$
The inequalities are equations whenever $v$ is in the support of $\mu$. We trivially have
\begin{align*}
 \| \phi(v) - \phi(w)\|_{\ell^1} &\leq \sum_{j=1}^{m} \mu(w_j)|d(v, w_j) -d(w,w_j)|\\
 &\leq \sum_{j=1}^m \mu(w_j) d(v,w) = d(v,w)
 \end{align*}
which shows that the embedding is $1-$Lipschitz with respect to $\ell^1(\mathbb{R}^m)$. Let us now fix two arbitrary vertices $w_1, w_2 \in \supp \mu$. The entry of $\phi(w_1) - \phi(w_2)$ in the $w_1-$th position is $-d(w_1, w_2) \mu(w_2)$. This implies
$$ \| \phi(w_1) - \phi(w_2) \|_{\ell^{\infty}(\mathbb{R}^m)} \geq d(w_1, w_2) \mu(w_2).$$
This inequality is trivially correct also when $w_1 = w_2$. Thus, summing $w_2$ over all vertices in the boundary, we get
$$ \frac{1}{m} \sum_{w_2 \in \supp(\mu)}  \| \phi(w_1) - \phi(w_2) \|_{\ell^{\infty}(\mathbb{R}^m)} \geq  \frac{1}{m} \sum_{w_2 \in \supp(\mu)} d(w_1, w_2) \mu(w_2) = \frac{\alpha}{m}.$$
Recalling that $\alpha \geq \diam(G)/2$, we obtain the desired result.
\end{proof}

\subsection{Proof of Proposition 1}
\begin{proof} We fix $D \in \mathbb{R}^{n \times n}$ to be the distance matrix. Throughout the proof we will identify probability measures with probability vectors. We first show that critical points of the functional are balanced. Suppose now $\mu$ is a critical point and $\nu$ is an arbitrary signed measure
with total weight $\nu(V) = 0$ such that
$$ \supp \nu \subseteq \supp \mu = \left\{v \in V: \mu(v) > 0 \right\}.$$
Then, for some $\varepsilon_0 > 0$ (depending only on $\mu$ and $\nu$) and all $|\varepsilon| \leq \varepsilon_0$, we have that $\mu + \varepsilon \nu$ is also a probability measure and thus, since $\mu$ is a critical point,
$$ \left\langle (\mu + \varepsilon \nu), D(\mu + \varepsilon \nu) \right\rangle = \left\langle \mu, D\mu \right\rangle + o(\varepsilon).$$
Since $D$ is symmetric, the linear term of the left-hand side is given by
$$ \left\langle \nu, D\mu \right\rangle +  \left\langle D\nu, \mu \right\rangle = 2\left\langle D\mu, \nu\right\rangle = 0.$$
This, in turn, implies that $D\mu$ restricted to $\supp \mu$ has to be constant since otherwise we could construct a signed measure $\nu$ for which the equation is not satisfied which shows that $\mu$ could not have been a critical point.  Suppose now that this constant $c$, the restriction of $D\mu$ onto $\supp \mu$, is different from $\|D\mu\|_{\ell^{\infty}}$ (in which case it has to satisfy $c < \|D \mu \|_{\ell^{\infty}}$. Then there exists a vertex $w \in V$ with $w \notin \supp \mu$ and $(D\mu)(w) > c$. Then, by taking any $ w_2 \in \supp(\mu)$, we can construct the measure $\nu = \delta_w - \delta_{w_2}$ and see that, for $\varepsilon$ sufficiently small, $\mu + \varepsilon \nu$ is still a probability measure with larger energy which again contradicts the fact that $\mu$ was a critical point. This shows that $D\mu$ equals $\|D \mu \|_{\ell^{\infty}}$ when restricted to $\supp \mu$ and is therefore a balanced measure. Let us now suppose that $\mu$ is a balanced measure. Let $\nu$ be an arbitrary signed measure such that $\mu + \varepsilon \nu$ is a probability measure for all $\varepsilon$ sufficiently small (this means that the negative entries of $\nu$ have to be contained in $\supp \mu$). Then
$$ \left\langle (\mu + \varepsilon \nu), D(\mu + \varepsilon \nu) \right\rangle = \left\langle \mu, D\mu \right\rangle + 2 \varepsilon \left\langle D\mu, \nu \right\rangle + o(\varepsilon).$$
Since the negative weights of $\nu$ are contained in $\supp \mu$ and since $D\mu$ is maximal in $\supp \mu$, we have that $ \left\langle D\mu, \nu \right\rangle \leq 0$ with equality if and only if all the positive weight of $\nu$ is also contained points where $D\mu$ assumes its maximum. This implies the statement. \end{proof}

\subsection{Proof of Proposition 2}
\begin{proof} Let $G=(V,E)$ be fixed, let $x_1, \dots, x_k \in V$ be given and suppose that $v \in V$ is \textit{not} a boundary vertex. We will then prove that
$$ \sum_{i=1}^{k} d(v,x_i) < \max_{w \in V}  \sum_{i=1}^{k} d(w,x_i)$$
which establishes the desired result. Since $v$ is not a boundary vertex, this means that for each vertex $x_i$, 
$$ \frac{1}{\deg(v)} \sum_{(v, w) \in E} d(w,x_i) \geq d(v,x_i).$$
Summing over all $x_i$, we get that
$$ f(v) = \sum_{i=1}^{k}  d(v,x_i) \leq \frac{1}{\deg(v)} \sum_{(v, w) \in E} \sum_{i=1}^{k} d(w,x_i) =   \frac{1}{\deg(v)} \sum_{(v, w) \in E} f(w).$$
This shows that the value of $f$ in a vertex $v \notin \partial G$ can be bounded from above by the average value in an adjacent vertex. Therefore,
if $f$ assumes its maximum in any non-boundary vertex, it also assumes the maximum in all adjacent vertices. It is then possible to hop from vertex
to vertex until one ends up being adjacent to a boundary vertex and thus the maximum is also assumed at the boundary.
\end{proof}

\end{document}